\documentclass{gtart}

\def\ifplaintex{\expandafter\ifx\csname documentclass\endcsname\relax}

\ifplaintex 
\else
\fi

\expandafter\ifx\csname beginpicture\endcsname\relax
\expandafter\ifx\csname documentclass\endcsname\relax
\input pictex \else
\input prepictex \input pictex \input postpictex \fi\fi

\def\gt{{\mathsurround=0pt\it $\cal G\mskip-2mu$eometry \&\ 
$\cal T\!\!$opology}}        

\def\gtp{{\mathsurround=0pt\it $\cal G\mskip-2mu$eometry \&\ 
$\cal T\!\!$opology $\cal P\!$ublications}}  

\def\lognumber#1{\def\thelognumber{#1}}
\def\volumenumber#1{\def\thevolumenumber{#1}}
\def\papernumber#1{\def\thepapernumber{#1}}
\def\volumeyear#1{\def\thevolumeyear{#1}}

\def\pagenumbers#1#2{\def\startpage{#1}\def\finishpage{#2}}
\def\published#1{\def\publishdate{#1}}
\def\proposed#1{\def\theproposer{#1}}
\def\seconded#1{\def\theseconders{#1}}
\def\received#1{\def\receiveddate{#1}}
\def\revised#1{\def\reviseddate{#1}}
\def\accepted#1{\def\accepteddate{#1}}
\def\asciititle#1{\def\theasciititle{#1}}

\long\def\asciiabstract#1{\long\def\theasciiabstract{#1}}

\let\\\par\let\thelognumber\relax
\let\thevolumenumber\relax\let\thepapernumber\relax
\let\thevolumeyear\relax\let\thesamplenumber\relax\let\startpage\relax
\let\finishpage\relax\let\publishdate\relax\let\receiveddate\relax
\let\reviseddate\relax\let\accepteddate\relax\let\theasciititle\relax
\let\theasciiauthors\relax
\let\theasciiabstract\relax
\let\theasciiemail\relax\let\theshortauthors\relax\let\theshorttitle\relax

\long\def\maketitlep{   

\count0=\startpage

\gt\hfill      
\beginpicture
\setcoordinatesystem units <0.33truein, 0.33truein> point at 2.2 0.9
\setplotsymbol ({$\cal G$})
\plotsymbolspacing=9truept
\circulararc 315 degrees from 0 1 center at 0 0
\setplotsymbol ({$\cal T$})
\circulararc 315 degrees from 1 -1 center at 1 0
\endpicture
\break
{\small\ifx\thesamplenumber\relax 
Volume \else Sample
\fi\thevolumenumber\ (\thevolumeyear)
\startpage--\finishpage\nl
Published: \publishdate}
\vglue 0.5truein plus 0.4fil minus 0.1truein

{\parskip=0pt\leftskip 0pt plus 1fil\def\\{\par\smallskip}{\ifplaintex\large
\else\Large\fi\bf\thetitle}\par\medskip}   

\vglue 0pt plus 0.1fil 

{\parskip=0pt\leftskip 0pt plus 1fil\def\\{\par}{\sc\theauthors}
\par\medskip}

\vglue 0pt plus 0.1fil 

{\small\parskip=0pt\let\newline\\
{\leftskip 0pt plus 1fil\def\\{\par}{\sl\theaddress}\par}
\expandafter\ifx\theemail\relax    
\relax\else\vglue 5pt plus 0.02fil minus 2pt\def\\{\stdspace{\rm 
and}\stdspace} 
\cl{Email:\stdspace\tt\theemail}\fi
\ifx\theurl\relax                  
\relax\else\vglue 5pt plus 0.02fil minus 2pt\def\\{\stdspace{\rm 
and}\stdspace}
\cl{URL:\stdspace\tt\theurl}\fi\par}

\vglue 7pt plus 0.3fil minus 3pt

{\bf Abstract}
\vglue 5pt plus 0.1fil minus 2pt

\theabstract

\vglue 7pt plus 0.3fil minus 3pt

{\bf AMS Classification numbers}\quad Primary:\quad \theprimaryclass

Secondary:\quad \thesecondaryclass

\vglue 5pt plus 0.3fil minus 2pt

{\bf Keywords}\quad \thekeywords

\vglue 10pt plus 0.5fil minus 5pt

{\small  Proposed: \theproposer\hfill Received: \receiveddate\nl
Seconded: \theseconders\hfill 
\ifx\reviseddate\relax                         
Accepted: \accepteddate                        
\else
Revised: \reviseddate                          
\fi}
\eject
}       

\let\maketitlepage\maketitlep
\let\maketitle\maketitlepage

\font\phead=cmsl9 scaled 950
\font=cmsl9 scaled 1050
\font\pnum=cmbx10 scaled 913
\font\lnum=cmbx10 
\font\pfoot=cmsl9 scaled 950
\font=cmsl9 scaled 1050
\ifplaintex
\headline{\vbox to 0pt{\vskip -4.5mm\line{\small\phead\ifnum
\count0=\startpage ISSN 1364-0380 (on line)
1465-3060 (printed) \hfill {\pnum\folio}\else\ifodd\count0\def\\{ }%
\ifx\theshorttitle\relax\thetitle\else\theshorttitle\fi\hfill{\pnum\folio}
\else\def\\{ and }{\pnum\folio}\hfill\ifx\theshortauthors\relax\theauthors
\else\theshortauthors\fi\fi\fi}\vss}}
\footline{\vbox to 0pt{\vglue 0mm\line{\small\pfoot\ifnum\count0=\startpage
\copyright\ \gtp\hfill\else
\gt, Volume \thevolumenumber\ (\thevolumeyear)\hfill\fi}\vss
}}
\else
\makeatletter
\def\@oddhead{{\small\lhead\ifnum\count0=\startpage ISSN 1364-0380 (on line)
1465-3060 (printed) \hfill {\lnum\number\count0}\else\ifodd\count0
\def\\{ }\ifx\theshorttitle\relax \thetitle \else\theshorttitle\fi\hfill
{\lnum\number\count0}\else\def\\{ and }{\lnum\number\count0}
\hfill\ifx\theshortauthors\relax 
\theauthors\else\theshortauthors\fi\fi\fi}}\def\@evenhead{\@oddhead}
\def\@oddfoot{\small\lfoot\ifnum\count0=\startpage\copyright\ \gtp\hfill\else
\gt, Volume \thevolumenumber\ (\thevolumeyear)\hfill\fi}
\def\@evenfoot{\@oddfoot}
\makeatother
\fi

\newwrite\gtoutfile
\long\gdef\makeheadfile{  
{\def\\{, }\def\s{ }
\immediate\openout\gtoutfile head.xxx
\immediate\write\gtoutfile{Proxy-for: \ifx\theasciiauthors\relax
\theauthors\else\theasciiauthors\fi\s<\ifx\theasciiemail\relax\theemail\else\theasciiemail\fi>}
\immediate\write\gtoutfile{\noexpand\\}
\immediate\write\gtoutfile{Authors: \ifx\theasciiauthors\relax
\theauthors\else\theasciiauthors\fi}
{\def\\{ }\immediate\write\gtoutfile{Title: \ifx\theasciititle\relax
\thetitle\else\theasciititle\fi}}
\immediate\write\gtoutfile{Subj-class: GT or SG or MG etc}
\immediate\write\gtoutfile{MSC-class: \theprimaryclass\ifx\thesecondaryclass\relax\else, \thesecondaryclass\fi}
\immediate\write\gtoutfile{Journal-ref: Geom. Topol. \thevolumenumber
(\thevolumeyear) \startpage-\finishpage}
\immediate\write\gtoutfile{Comments: Published by Geometry and Topology at}
\immediate\write\gtoutfile{\s\s http://www.maths.warwick.ac.uk/gt/GTVol\thevolumenumber/paper\thepapernumber.abs.html}
\immediate\write\gtoutfile{\noexpand\\}
\immediate\write\gtoutfile{}
\ifx\theasciiabstract\relax
\immediate\write\gtoutfile{\theabstract}\else
\immediate\write\gtoutfile{\theasciiabstract}\fi
\immediate\write\gtoutfile{}
\immediate\write\gtoutfile{\noexpand\\}
\immediate\write\gtoutfile{}
\immediate\closeout\gtoutfile}}  

\def\maketitlepage{\maketitlep\makeheadfile}
\let\maketitle\maketitlepage

\lognumber{299}
\volumenumber{7}\papernumber{28}\volumeyear{2003}
\pagenumbers{965}{999}
\received{9 January 2003}
\revised{10 December 2003}
\published{21 December 2003}
\accepted{19 December 2003}
\proposed{Ronald Stern}
\seconded{Ronald Fintushel, John Morgan}

\def\S{section~}
\let\oldcolon\colon
\def\colon{\oldcolon\thinspace}
\makeatletter

\@addtoreset{equation}{section}
\makeatother

\renewcommand{\theenumi}{\roman{enumi}}

\newtheorem{theorem}{Theorem}[section]
\newtheorem{prop}[theorem]{Proposition}
\newtheorem{lemma}[theorem]{Lemma}

\newtheorem{cor}[theorem]{Corollary}
\newtheorem{remark}[theorem]{Remark}
\newtheorem{defn}[theorem]{Definition}

\newcommand{\R}{{\bf R}}
\newcommand{\C}{{\bf C}}

\newcommand{\Z}{{\bf Z}}
\newcommand{\qu}{{\bf H}}
\newcommand{\A}{{\cal A}}
\newcommand{\G}{{\cal G}}
\newcommand{\B}{{\cal B}}

\newcommand{\K}{{\cal K}}
\newcommand{\M}{{\cal M}}
\newcommand{\N}{{\cal N}}

\newcommand{\CC}{{\cal C}}
\newcommand{\I}{{\rm I}}

\newcommand{\XX}{{\cal X}}

\newcommand{\stab}{\mbox{\rm stab}}
\newcommand{\kernel}[1]{{\rm ker\,}{#1}}

\newcommand{\proofof}[1]{\proof[Proof of #1]}

\newcommand{\tr}{\mbox{\rm Tr}}

\newcommand{\normi}[2]{\|#1\|_{#2}}
\newcommand{\normii}[3]{\|#1\|_{#2}^{#3}}
\newcommand{\const}{\mbox{\rm const.}}

\newcommand{\dimn}{{\rm dim}}

\newcommand{\cs}{{\rm cs}}
\newcommand{\csd}{{\rm csd}}
\newcommand{\su}{\mbox{\rm su}}
\newcommand{\SF}{{\rm SF}}

\newcommand{\dis}{\displaystyle}
\newcommand{\supp}{{\rm supp\,}}
\newcommand{\ad}{{\rm ad}}
\newcommand{\Ad}{{\rm Ad}}
\newcommand{\detind}{{\rm detind\,}}

\renewcommand{\epsilon}{\varepsilon}
\newcommand{\mPhi}{{\mit\Phi}}
\newcommand{\mPsi}{{\mit\Psi}}
\newcommand{\mLambda}{{\mit\Lambda}}

\newcommand{\mOmega}{{\mit\Omega}}
\newcommand{\mTheta}{{\mit\Theta}}

\begin{document}

\title[A non-abelian Seiberg--Witten invariant]{A non-abelian 
Seiberg--Witten invariant for\\integral homology 3--spheres}
\author{Yuhan Lim}
\email{ylim@math.ucsc.edu}
\address{Dept of Mathematics,
Univ of California\\Santa Cruz, CA 95064, USA}

\asciititle{A non-abelian Seiberg-Witten invariant for 
integral homology 3-spheres}
\begin{abstract}
A new diffeomorphism invariant of integral homology 
3--spheres is defined using a non-abelian ``quaternionic''
version of the Seiberg--Witten equations.
\end{abstract}

\asciiabstract{A new diffeomorphism invariant of integral homology 
3-spheres is defined using a non-abelian 'quaternionic'
version of the Seiberg-Witten equations.}




\maketitlepage

\section{Introduction}

The Seiberg--Witten equations when applied to the study of
oriented integral homology 3--spheres yield an invariant
which was shown in \cite{lim1} to coincide with Casson's
invariant. In \cite{BH}, Boden and Herald introduced a 
generalization of Casson's invariant from $SU(2)$ to the higher structure
group $SU(3)$ based on the gauge theory approach of Taubes \cite{taubes}. 
This $SU(3)$--Casson invariant 
utilizes values of the Chern--Simons function which 
makes it a real valued invariant rather than an integral one.
In the present article we define a non-Abelian
version of the Seiberg--Witten equations which we call
\textit{quaternionic} and construct a topological invariant 
of integral homology 3--spheres in a  
manner parallel to the $SU(3)$--Casson invariant.
This new invariant has the property that it is independent of 
orientation of the 3--manifold and
a linear combination with the $SU(3)$--Casson invariant
gives a $\Z\bmod 4\Z$ invariant for unoriented integral homology 3--spheres. 

The contents of this article are as follows.
In \S\ref{gauge setup} we introduce the generalization of the
SW--equations we use. The technical issue of admissible perturbations is
also discussed. We use the novel approach of non-gradient perturbations.
Section \ref{spectral} gives the main results which are
Theorems~\ref{thrm1} and \ref{thrm2}.
The remaining sections take up the proofs.
We assume the reader has some familiarity
with  \cite{BH}, \cite{lim1} and \cite{taubes}.

\section{Quaternionic gauge theory in 3--dimensions}\label{gauge-setup}\label{gauge
setup}

\textbf{Standing Convention} \textit{Throughout this article {\it\,Y} will denote an oriented closed 
integral homology $3$--sphere (ZHS). {\it\,Y} will also be assumed to
have a fixed Riemannian metric $g$.} 

The aim is to introduce a quaternionic setting in which the
Seiberg--Witten equations will make sense. 
Since $Y$ is a ZHS it has a unique
spin structure, up to equivalence. With respect to $g$ this is
given by a principal $\mbox{\sl spin}(3)\cong \mbox{\sl SU}(2)$ bundle $P\to Y$. In the (real) Clifford
bundle $CL(T^{*}Y)\cong CL(Y)$ the volume form $\omega_Y$ has the property
that $\omega^2_{Y}=1$. The action of $\omega_Y$ on $CL(Y)$ induces
a splitting into $\pm 1$ eigenbundles $CL^{+}\oplus CL^{-}$.
Both $CL^{+}$ and $CL^{-}$ are bundles of algebras over $Y$
with each fibre isomorphic, as an algebra, to the quaternions $\qu$.

Let $S\to Y$ be the complex spinor bundle on which $CL^{+}$ acts non-trivially.
This is a rank 2 complex Hermitian vector bundle. Since the
fibres of $CL^{+}$ are quaternionic vector spaces,
$S$ possesses an additional action by $\qu$ which commutes with the
Clifford action (see \cite{Law}); we may take this to be a right
action $S\times{\qu}\to S$.

Suppose now that $E\to Y$ is a given fixed rank one metric quaternionic
vector bundle --- we assume the action by $\qu$ is a left action. $E$ also has a description as
a complex Hermitian rank 2 vector with trivial determinant,
i.e.\ with structure group  $SU(2)$.
We can twist the spinor bundle $S$ by tensoring with $E$ over
the quaternions to form the bundle $S\otimes_{\qu}E$.
This is a {\it real} rank 4 Riemannian vector bundle and does not
naturally inherit a complex structure from $S$ or $E$. 

Given an $SU(2)$--connection $A$ on $E$ (henceforth any connection
on $E$ mentioned will be assumed to be such type) we may construct
using the canonical Riemiannian connection on $S$, a metric
(i.e.\ $SO(4)$) connection on $S\otimes_{\qu}E$. This then
defines in the usual way a Dirac operator
$$
D_A = \sum_{i=1}^3 e_i\cdot \nabla^{A}_{e_i}.
$$
Here the $e_i$ are an orthonormal frame and $\nabla^A$ is the connection
on $S\otimes_{\qu}E$ mentioned above. We
emphasize that $D_A$ is in general only a real linear operator on
$S\otimes_{\qu}E$.

\begin{lemma}\label{complexification}
The complexification of $S\otimes_{\qu}E$ is naturally isomorphic
as a complex Clifford module with $S\otimes_{\C}E$.
Under this isomorphism the complexification $D_{A}\otimes\C$ corresponds to 
the complex Dirac operator $D^{\C}_{A}$.
\end{lemma}

\proof
Introduce the notation $\overline{\otimes}$ to denote the tensor product
of elements in $S\otimes_{\qu}E$ and $\otimes$ the complex tensor product
in $S\otimes_{\C}E$.
Define the vector bundle map $h$
from $S\otimes_{\C}E$ to 
$(S\otimes_{\qu}E)\otimes\C$ by 
$$
h(e\otimes f)  
=e\overline{\otimes} f-\sqrt{-1}(e i\overline{\otimes} f).
$$
One checks directly that this map is a complex isomorphism and commutes with
Clifford multiplication. \endproof

Since the real two forms ${\mLambda}^2$ naturally include in $CL(Y)$ we have 
by Clifford mutiplication the action of ${\mLambda}^2$ on $S$.
This representation of ${\mLambda}^2$ on $S$
is well-known to be injective and with image the adjoint bundle $\ad S$, the bundle of
skew-Hermitian transformations of $S$.
The bundle $\ad E$ acts on $E$ from the left.
Define an action of
${\mLambda}^2\otimes \ad E$ on $S\otimes_{\qu}E$ by the
rule
$$
(\omega\otimes l)\cdot (\phi\otimes e) := 
(\omega\cdot\phi)\otimes l(e).
$$
This is well-defined since the actions of $\mLambda^2$ and $ \ad E$
commute with the quaternionic structures.

\begin{remark}\rm
The Clifford action of $\beta\in\mLambda^{2}$ is the same as the
action of $-*\beta\in\mLambda^{1}$ on $S$ since the volume form $\omega_{Y}$
acts by the identity. Thus we may equivalently work (up to multiplication
by $-1$) with the action of $\mLambda^{1}\otimes\ad E$ on $S\otimes_{\qu}E$.
\end{remark}

\begin{lemma}
The representation ${\mLambda}^2\otimes \ad E\to
\mbox{\rm End}_{\R}(S\otimes_{\qu}E)$ above is injective and
has image the subbundle $\mbox{\rm Sym}^0_{\R}(S\otimes_{\qu}E)$ 
of trace zero real symmetric transformations of $S\otimes_{\qu}E$.
\end{lemma}

\proof  
That the representation of the lemma is injective is
easily verified. We may rewrite the action of $\mLambda^2\otimes \ad E$
as $(\omega\otimes l)\cdot(\phi\otimes e)=- (i\omega\cdot \phi)\otimes il(e)$.
Since $i \ad S$ is exactly the trace zero Hermitian
symmetric bundle endomorphisms of $S$, and similiarly for
$i  \ad E$, the image of the representation clearly
lies in the trace zero real symmetric endomorphisms of
$S\otimes_{\qu}E$. That it is onto follows by a dimension
count giving both ${\mLambda}^2\otimes \ad E$ and
$\mbox{\rm Sym}^0_{\R}(S\otimes_{\qu}E)$ real vector bundles of rank $9$.
\endproof

The above lemma shows that we may regard the bundle 
$\mbox{\rm Sym}^0_{\R}(S\otimes_{\qu}E)$ as identical to
${\mLambda}^2(Y)\otimes \ad E$. Thus whenever convenient
we can think of a trace zero real symmetric endomorphism of
$S\otimes_{\qu}E$ as a twisted 2--form with values in $ \ad E$.

\begin{lemma}\label{bilinear}
There is a unique fibrewise symmetric bilinear form $\{\cdot\}_0$ on
$S\otimes_{\qu}E$ with values in ${\mLambda}^2\otimes \ad E$
determined by the rule that
$$
\langle{\omega,\{\phi\cdot\psi\}_0 }\rangle
=\langle{\omega\cdot\psi,\phi}\rangle 
=\langle{\omega\cdot\phi,\psi}\rangle
$$
holds for all sections $\omega$ of ${\mLambda}^2\otimes \ad E$.
As a section of $\mbox{\rm Sym}^0_{\R}(S\otimes_{\qu}E)$,
$\{\phi\cdot\psi\}_0$ is given by the expression
$$
\{\phi\cdot\psi\}_0=\frac{1}{2}\left(\phi\otimes\psi^*+\psi\otimes\phi^*
-\frac{1}{2}\langle{\phi,\psi}\rangle \I\right).
$$
Here $\phi\otimes\psi^*(\nu)=\phi\langle{\nu,\psi}\rangle$ and similiarly
for $\psi\otimes\phi^*$.
\end{lemma}

\proof Let $\{\phi_{i}\}$, $\{\omega_{j}\}$  be a local  orthonormal frames for 
$S\otimes_{\qu}E$, $\mLambda^{2}\otimes\ad E$ respectively. Let 
$\{\phi_{i}\cdot\phi_{j}\}_{0}=c_{i,j}^{k}\omega_{k}$. Then we see 
that $c_{i,j}^{k}=\langle{\omega_{k}\cdot\psi_{i},\psi_{j}}\rangle 
=\langle{\omega_{k}\cdot\psi_{j},\psi_{i}}\rangle=c_{j,i}^{k}$ determines 
$\{\cdot\}_{0}$.
Identify $\mLambda^{2}\otimes\ad E$ with $\mbox{\rm Sym}^0_{\R}(S\otimes_{\qu}E)$. 
In a local trivialization we may regard 
sections of $\mbox{\rm Sym}^0_{\R}(S\otimes_{\qu}E)$ as functions with values in $\mbox{\rm Sym}^0_{\R}(\R^{4})$, 
the $4\times 4$ real symmetric matrices, and sections of 
$S\otimes_{\qu}E$ as $\R^{4}$--valued functions.
As such the inner product in $\mbox{\rm Sym}^0_{\R}(\R^{4})$ is given by $\langle M,N\rangle = \tr(MN)$.
The right side of the defining equation for $\{\cdot\}_{0}$ can be 
expressed locally as 
\begin{eqnarray*}
\frac{1}{2}\left(\tr(W\mPhi\mPsi^{T})+\tr(W\mPsi\mPhi^{T})\right)
&=&\frac{1}{2}\tr(W(\mPhi\mPsi^{T}+\mPsi\mPhi^{T}))\\
&=&\langle{W,\frac{1}{2}(\mPhi\mPsi^{T}+\mPsi\mPhi^{T})}\rangle.
\end{eqnarray*}
The claimed expression for $\{\cdot\}_{0}$ is exactly the trace-free 
component of the symmetric expression $\frac{1}{2}(\mPhi\mPsi^{T}+\mPsi\mPhi^{T})$. \endproof

The {\it configuration space} $\CC$ is the space of all pairs
$(A,\mPhi)$ consisting of an $SU(2)$--connection $A$ on $E$ and
a section $\mPhi$ of $S\otimes_{\qu}E$. As usual we should
work within the framework of a certain functional space; for us
choose $A$ and $\mPhi$ to be of class $L^2_2$ (for $A$ this means
$A-A_0$ is $L^2_2$ where $A_0$ is a fixed $C^{\infty}$--connection).
$\CC$ is  an affine space modelled on the Hilbert space
$$
L^2_2({\mLambda}^1\otimes\ad E)\times L^2_{2}(S\otimes_{\qu}E).
$$
The {\it gauge automorphism group} $\G$ in this case will consist of the $L^2_3$--bundle
automorphisms which preserve the quaternionic structure of $E$,
or equivalently the $L^2_3$--sections of $\Ad E$. Since $L^2_2\subset C^0$ in dimension 3,
$\CC$ and $\G$ consists of continuous objects.
$\G$ acts on $\CC$ by
$g\cdot(A,\mPhi)=(g(A),g^{-1}\mPhi)$. This action is differentiable
and the quotient we denote by $\B$. Our convention is that $g(A)$
is the pull-back of $A$ by $g$.

We have the following  observation: the stabilizer
\begin{displaymath}
\stab(A,\mPhi)=\left\{
\begin{array}{ccl}
\{1\} & \rm if &\mPhi\neq 0\\
\stab(A) &\rm if &\mPhi=0.
\end{array}\right.
\end{displaymath}
The possible choices for $\stab(A)$ are $\{\pm 1\}$, $U(1)$ or $SU(2)$. Note that in the
last possibility $A$ is necessarily a trivial  connection.
The pair $(A,\mPhi)$ is {\it irreducible} if $\mPhi\neq 0$ and {\it reducible}
otherwise. Thus $\G$ acts freely on $\CC^*$, the irreducible portion of $\CC$ and the
quotient $\CC^*$ by $\G$ is denoted $\B^*$.

$\G$ is a Hilbert Lie group
with tangent space at the identity
$
T_{e}\G=L^2_3(\ad E).
$
Let $\G\to\CC$, $g\mapsto (g(A),g^{-1}\mPhi)$ be the map which is the orbit 
of $(A,{\mPhi})$
under the action of $\G$. The derivative at the
identity is the map
\begin{equation}
\begin{array}{c}
\delta_{A,{\mPhi}}^0\colon L^2_3(\ad E)\to 
L^2_2({\mLambda}^1\otimes  \ad E)\oplus
L^2_2(S\otimes_{\qu}E),\\
\delta^0_{A,\mPhi}(\gamma) =(d_A\gamma,-\gamma(\mPhi) ).
\end{array}
\end{equation}
A slice for the action of $\G$ on $\CC$ at $(A,\mPhi)$
is given by $(A,\mPhi)+X_{A,\mPhi}$ where $X_{A,\mPhi}$ is the
{\it slice space} which is the $L^2$--orthogonal complement 
in $L^2_2({\mLambda}^1\otimes \ad E)\oplus L^{2}_{2}(S\otimes_{\qu}E)$ 
of the image of $\delta^{0}_{A,\mPhi}$. We may also regard $X_{A,\mPhi}$
as the tangent space to $\B^{*}$ at an irreducible orbit $[A,\mPhi]$.

Define a bilinear product 
$B\colon(S\otimes_{\qu}E)\otimes(S\otimes_{\qu}E)\to \ad E$ by the rule
that $\langle{\gamma(\phi),\psi}\rangle=\langle{\gamma,B(\psi,\phi)}\rangle$
holds for all $\gamma\in \ad E$. Then $X_{A,\mPhi}$ has the
description as the subspace of $L^2_2(\mLambda^1\otimes \ad E)\oplus
L^2_2(S\otimes_{\qu}E)$ defined by the equation
\begin{equation}\label{slice-equ}
\delta^{0*}_{A,\mPhi}(a,\psi)=0\quad \Longleftrightarrow\quad
d^*_A a - B(\mPhi,\psi)=0.
\end{equation}

A reducible we will often simply denote by
$A$ instead of $(A,0)$. Corresponding reducible subspaces of $\CC$ and $\B$
are denoted $\A$ and $\B_{\A}$.
At a reducible $A$ the slice $X_{A}$ splits into a product 
$X^{r}_{A}\times L^2_2(S\otimes_{\qu}E)$ where $X^{r}_A$ is  the slice
for the action of $\G$ on $\A$. Then the normal space  to
$\B_{\A}$ in $\B$ near $[A]$ is modelled on
$$
L^2_2(S\otimes E)/\stab(A).
$$
For instance if $A$ is irreducible as a connection then this normal
space is a cone on the quotient of the unit sphere in a separable Hilbert space
by the antipodal map $v\mapsto -v$.

On $\CC$ we have the Chern--Simons--Dirac function $\csd\colon\CC\to 
\R$ 
(with respect to a choice of trivial connection $\mTheta$ say) given by
$$
\csd(A,\mPhi) =\frac{1}{8\pi^2}\int_Y {\rm Tr}\Bigl( a\wedge d_{\mTheta}a 
+\frac{2}{3}a \wedge a\wedge a\Bigr)
-\int_Y \langle{D_A\mPhi,\mPhi}\rangle, \quad a=A-\mTheta.
$$
A direct computation gives
\begin{eqnarray*}
d\csd_{A,\mPhi}(a,\phi) &=& \int_Y {\rm Tr}(F_A\wedge a) 
+\int_Y \langle{*a\cdot\mPhi,\mPhi}\rangle
-\int_Y\langle{D_A\mPhi,\phi}\rangle\\
&=& \int_Y\langle{-*F_A+*\{\mPhi\cdot\mPhi\}_0,a}\rangle
-\int_Y\langle{D_A\mPhi,\phi}\rangle.
\end{eqnarray*}
Thus the negative of the $L^2$--gradient of $\csd$ is the \lq $L^2_{1}$--vector field\rq\ on $\CC$
\begin{equation}\label{chi-map}
{\cal X}(A,\mPhi) \stackrel{\rm def}{=} (*F_A -*\{\mPhi\cdot\mPhi\}_0, D_A\mPhi)
\in L^{2}_{1}.
\end{equation}
By this we mean that ${\cal X}$ is a section of the $L^2_{1}$--version of
the tangent bundle to $\CC$. The {\it Quaternionic Seiberg--Witten
equation} is the equation for the zeros of $\cal X$, i.e.\ the
critical points of $\csd$.

\begin{defn}\rm  
The {\it Quaternionic Seiberg--Witten equation} is
the equation defined for a pair $(A,\mPhi)$ consisting of a connection
on $E$ and a section $\mPhi$ (\lq{\it spinor}\rq) of
$S\otimes_{\qu}E$. The equation reads:
\begin{equation}\label{sw-equ}
\left\{
\begin{array}{ccl}
\dis F_A -\{\mPhi\cdot\mPhi\}_0 &=&0\\
D_A \mPhi &=& 0
\end{array}
\right.
\end{equation}
where $F_A$ is the curvature of $A$, and since $A$ is an 
$SU(2)$--connection, a section of ${\mLambda}^2\otimes  \ad E$.
$D_{A}$ is the Dirac operator on $S\otimes_{\qu}E$ and $\{\ \}_{0}$
denotes the quadratic form of Lemma~\ref{bilinear}.
\end{defn}

If $g$ is gauge transformation  then $\csd (g(A),g^{-1}\mPhi) = 
\csd(g,\mPhi) \pm {\rm deg}(g)$, so $\csd$ descends to an $\R/\Z$--valued
function on $\B$. This implies that 
${\cal X}(A,\mPhi) \in X_{A,\mPhi}\cap L^2_{1}$
and the portion of $\cal X$ over $\CC^{*}$ descends to a \lq $L^{2}_{1}$--vector
field\rq\ $\widehat{\cal X}$ over $\B^{*}$.

\begin{defn}\rm
The {\it moduli space} of solutions to (\ref{sw-equ}) we denote by
\begin{displaymath}
\M\stackrel{\rm def}{=}\{(A,\mPhi)\ {\rm solving}\ ({\rm \ref{sw-equ}})\}/\G\subset\B.
\end{displaymath}
$\M^*$ will denote irreducible  and $\M^{r}$ will denote the reducible portion of $\M$
respectively.
\end{defn}

Thus $\M^{*}$ is the zeros of $\widehat{\cal X}$ and following Taubes, will be
the basis for defining a Poincare--Hopf index for $\B^{*}$. 

\begin{remark}\rm
In our Quaternionic SW--theory the
reducible portion $\M^{r}$ of $\M$ is just
the moduli space of flat $SU(2)$--connections on $Y$. This is the 
space dealt with by Taubes \cite{taubes} in the gauge theory approach
to Casson's invariant. 
\end{remark}

We need to now address the issue of an admissible class of perturbations which
will make $\M$ a finite number of \textit{non-degenerate points} (made precise
below) to apply the idea of a Poincare--Hopf index. Unlike the holonomy
perturbations used by Taubes and Boden--Herald which are gradient perturbations
we elect to perturb $\cal X$ directly rather than $\csd$; i.e.\ at the level of
vector fields, for this avoids a number of
technical problems which the author has presently no satisfactory
solution. This approach will be adequate for defining a Poincare--Hopf
index but not a Floer type homology theory where gradient perturbations
are required.

\begin{defn}\rm 
An {\it admissible perturbation} $\pi$ consists of a 
differentiable $\G$--equivariant map of the form
$(*k,l)\colon\CC\to L^2_2(\mLambda^1\otimes \ad E)\times L^2_2(S\otimes_{\qu} E)$ where
\begin{enumerate}
\item
$\pi_{A,\mPhi}=(*k_{A,\mPhi},l_{A,\mPhi})\in X_{A,\mPhi}$
\item
the linearization of $(*k,l)$ at $(A,\mPhi)$ is a bounded linear
operator $$(L\pi)_{A,\mPhi}\colon L^{2}_{2}(\mLambda^1\otimes \ad E)\oplus L^2_2(S\otimes_{\qu} E)
\to L^{2}_{2}(\mLambda^1\otimes \ad E)\oplus L^2_2(S\otimes_{\qu} E)$$
\item
there is a uniform bound 
$$
\normi{\pi_{A,\mPhi}}{L^{2}_{2,A}} \stackrel{\rm def}{=}
\sum_{i=0}^{2}\normi{(\nabla^{A})^{i}\pi_{A,\mPhi}}{L^{2}} \le C.
$$
\end{enumerate}
\end{defn}

\begin{remark}\rm 
In the unperturbed case, $\M$ can be easily shown to be compact.
The preceding uniform $L^2_{2}$--type bound requirement on the perturbation
is crucial to retain compactness of the moduli space for the perturbed equation below. This
is a gauge invariant bound.
\end{remark}

\begin{defn}\rm
The {\it perturbed Quaternionic Seiberg--Witten equations} are the 
equations
\begin{displaymath}\label{sw-eq2}
\left\{
\begin{array}{ccl}
\dis F_A -\{\mPhi\cdot\mPhi\}_0
+k_{A,\mPhi}&=&0\\
\dis D_A {\mPhi} +l_{A,\mPhi} &=& 0.
\end{array}
\right.
\end{displaymath}
The corresponding moduli space is denoted $\M_{\pi}$, the
irreducible portion $\M^{*}_{\pi}$
and the reducible portion $\M^{r}_{\overline{\pi}}$ where
$\overline{\pi}$ is the restriction to $\A$ or equivalently the
$k$--component of $\pi$. Note that when $\mPhi=0$,
$\stab(A)$--invariance forces $l_{A,0}=0$ and the only effective
portion of $\pi$ on $\A$ is the $k$--component.
\end{defn}

Let ${\cal X}_{\pi}={\cal X}+\pi$, the perturbation of ${\cal X}$. The linearization
at a zero $(A,\mPhi)$ is a map
\begin{eqnarray}
&(L\XX_{\pi})_{A,\mPhi}\colon L^{2}_{2}(\mLambda^1\otimes \ad E)\oplus L^2_2(S\otimes_{\qu} E)
\to X_{A,\mPhi}\cap L^{2}_{1} \label{chi-lin}\\
&(L\XX_{\pi})_{A,\mPhi}(a,\phi) =
(*d_{A}a -*\{\phi\cdot\mPhi\}_{0}, D_{A}\phi+a\cdot\mPhi) +(L\pi)_{A,\mPhi}(a,\phi).\nonumber
\end{eqnarray}

\begin{defn}\rm
Call $(A,\mPhi)$ or $[A,\mPhi]$ {\it non-degenerate} if $L\XX_{\pi}$
is surjective at $(A,\mPhi)$. $\M_{\pi}$ is {\it non-degenerate} if it consists entirely
of non-degenerate points. In this instance we also call $\pi$ {\it non-degenerate}.
The standard Kuranishi local model argument shows that a non-degenerate point
is isolated in $\B$. (This includes reducible points.)
\end{defn}

Fix a connection $\nabla^{0}$ and let $L^{2}_{2}$ denote the Sobolev
norm with respect to $\nabla^{0}$. A metric on $\B$ is defined by the
rule
\begin{equation}\label{metric-def}
d([A,\mPhi],[A',\mPhi']) =\inf_{g\in\G}\Bigl\{ \normi{(A-g(A'),\mPhi-g^{-1}\mPhi')}{L^{2}_{2}}\Bigr\}.
\end{equation}

\begin{prop}\label{adm-lem}
For any admissible perturbation $\M_{\pi}$ is a compact subspace
of $\B$. Furthermore there is an $\epsilon_{0}>0$ such that for any 
$0<\epsilon<\epsilon_{0}$,
if $\normi{\pi_{A,\mPhi}}{L^{2}_{2,A}}<\epsilon$ uniformly
then given any $[A,\mPhi]\in \M_{\pi}$ there is a $[A',\mPhi']\in\M$ such
that $d([A,\mPhi],[A',\mPhi'])<\epsilon$.
\end{prop}

\begin{prop}\label{adm2-lem}
There exists non-degenerate admissible perturbations.\break
Furthermore such a perturbation $\pi$ may be choosen so that
$\normi{\pi_{A,\mPhi}}{L^{2}_{2,A}}$
is arbitrarily small (uniformly)  and $\pi$
vanishes on any given closed subset of $\CC$ which is disjoint
from the subspace of unperturbed SW--solutions.
\end{prop}

The proofs are in sections \ref{comp-sect} and \ref{pert-sect}.

\section{Spectral flow and definition of the invariant}\label{spectral}

Fix a perturbation $\pi$ (not necessarily non-degenerate). 
Regard the image of $\XX_{\pi}$ as lying in the larger space
$L^2(\mLambda^1\otimes \ad E)\oplus L^2(S\otimes_{\qu} E)$
since $X_{A,\mPhi}\cap L^{2}_{1}$ is a subspace of the former.
The analog of the operator used by Taubes to define relative
signs between non-degenerate zeros of $\widehat{\XX_{\pi}}$
is the unbounded operator on 
$L^{2}( \ad E)\oplus L^2(\mLambda^{1}\otimes \ad E)\oplus L^2(S\otimes_{\qu} E)$
given in block matrix form:
\begin{displaymath}
L^{\pi}_{A,\mPhi} =
\left(
\begin{array}{cc}
0& \delta^{0 *}_{A,\mPhi} \\
\delta^{0}_{A,\mPhi}& (L\XX_{\pi})_{A,\mPhi}
\end{array}
\right).
\end{displaymath}
Here $\delta^{0 *}_{A,\mPhi}$ is the formal $L^{2}$--adjoint of $\delta^{0}_{A,\mPhi}$ 
and the splitting used above is the 
$(\mbox{\rm first})\oplus(\mbox{\rm 2nd and 3rd factors})$.
$L^{\pi}_{A,\mPhi}$ has dense domain the subspace of $L^{2}_{2}$--sections.
The ellipticity of $L^{\pi}_{A,\mPhi}$ implies that it is  closed and
unbounded as an operator on $L^{2}$. 
In general $L^{\pi}_{A,\mPhi}$
will not have a real spectrum, due to the non-gradient
perturbations we are using.

\begin{remark}\rm
If it were the case that $L^{\pi}_{A,\mPhi}$ is 
formally self-adjoint (i.e.\ on smooth sections)  
then it is  
well-known that $L^{\pi}_{A,\mPhi}$ has only a
discrete real spectrum which is unbounded in both directions in
$\R$ and is without any accumulation points.
It can be shown in general that since 
$L^{\pi}_{A,\mPhi}$ is a $L_{A,\mPhi}$--compact perturbation of
$L_{A,\mPhi}$ on $L^{2}$ the spectrum continues to be discrete,
the real part of the spectrum 
is also unbounded in both directions in $\R$ and is without any accumulation
points, see \cite{kato}.
\end{remark}

Let us consider the behaviour of $L^{\pi}_{A,\mPhi}$ along the
reducible stratum $\A\subset\CC$. Since $\mPhi=0$ we abbreviate
the operator to $L^{\pi}_{A}$. This has a natural splitting
\begin{equation}\label{split-eq}
L^{\pi}_{A} = K^{\overline{\pi}}_{A}\oplus D^{\pi}_{A}
\end{equation}
corresponding to the splitting $(\mLambda^{0+1}\otimes \ad E)\oplus (S\otimes_{\qu}E)$.
We call $K^{\overline{\pi}}_{A}$ the {\it tangential operator} (the dependence
on only the restriction $\overline{\pi}$ of $\pi$ will be clear below)
and $D^{\pi}_{A}$ the {\it normal operator}. Explicitly 
\begin{displaymath}
K^{\overline{\pi}}_{A} =
\left(
\begin{array}{cc}
0 & d^{*}_{A}\\ d_{A} & *d_{A}+(L\overline{\pi})_{A}
\end{array}
\right)
\quad\mbox{\rm and}\quad
D^{\pi}_{A} = D_{A}+(L^{\nu}\pi)_{A}.
\end{displaymath}
Here $(L^{\nu}\pi)_{A}$ denotes $(L\pi)_{A,0}$ restricted
to the normal space
$L^{2}_{2}(S\otimes_{\qu}E)$  followed by
projection onto the normal space again. The normal space is naturally
acted on by $\stab(A)$ and $D^{\pi}_{A}$ commutes with this action.
If  $A$ is irreducible as a connection then $\stab(A)=\{\pm 1\}$ and
$D^{\pi}_{A}$ remains real linear but when $A=\mTheta$ a trivial connection,
$\stab(A)\cong SU(2)$ and it is quaternionic linear. (There is 
a $\stab(A)\cong U(1)$ case but this will not play a role so we will
omit discussing it.)

A fact established in \S\ref{pert-sect} is:

\begin{lemma}\label{sf-fact-lem}
The operator $(L^{\nu}\pi)_{A}$ is multiplication by a real function $f_{A}\in L^{2}_{2}$.
Thus $D^{\pi}_{A}=D_{A}+f_{A}$ extends to an unbounded
self-adjoint operator on
$L^{2}(S\otimes_{\qu}E)$ and spectral flow is defined for this
operator.
\end{lemma}

To define relative signs between non-degenerate
zeros of $\widehat{\XX_{\pi}}$ one usually uses the $\bmod 2$
spectral flow of $L^{\pi}_{A,\mPhi}$ when this operator is
self-adjoint. In the general case we use the determinant
line $\detind L^{\pi}$ regarding $L^{\pi}_{A,\mPhi}$ as a family
parameterized by $(A,\mPhi)\in\CC$. This is equivalent to the spectral
flow definition in the self-adjoint case. $\detind L^{\pi}$
descends to a line bundle over $\B$ which we also denote by the
same notation. However we note:

\begin{lemma}
$\detind L^{\pi}$ is non-orientable over $\B$, i.e.\ there
exists closed loops $\gamma\colon S^{1}\to \B$ such that
$\gamma^{*}(\detind L^{\pi})$ is a non-trivial line bundle
over $S^{1}$.
\end{lemma}

\proof It suffices to consider the determinant index $\detind L$ of the unperturbed 
family over $\B_{\A}\subset\B$. Then $L_{A}=K_{A}\oplus D_{A}$ where $K_{A}$ is
essentially the boundary of the (twisted) Self-dual operator in dimension 4. 
Spectral flow around closed loops for $K_{A}$ is equivalent to the index 
of the (twisted) Self-dual operator over $Y\times S^{1}$. The latter index is well-known
to be $\equiv 0\bmod 8$. By Lemma~\ref{complexification} the spectral
flow for $D_{A}$ around closed loops is equivalent to the index of the twisted complex
Dirac operator over $Y\times S^{1}$, where the twisting bundle $E$ is rank 2 complex.
According to the Atiyah--Singer Index Theorem this index is the negative of 2nd Chern class
of $E$ evaluated over the fundamental class of $Y\times S^{1}$. 
We may choose any closed loop so that this is $\pm 1$. (See the proof 
of Lemma~\ref{sf-calc}
for more details on this part of the calculation.)\endproof

In particular if we have two
non-degenerate zeros of $\widehat{\XX}_{\pi}$ then the Lemma asserts
that it is impossible in this scheme to define a relative i.e.\ $\bmod 2$ 
sign between non-degenerate zeros.  Thus as far as defining an invariant goes
we can only work with the cardinality 
$$
\sum_{[A]\in \M_{\pi}^{*}} 1 \quad \bmod 2
$$
for $\pi$ non-degenerate.

Assume now that $\pi$ is non-degenerate. We define counter-terms
associated to $\M^{r}_{\overline{\pi}}$
to make $\sum_{[A]\in \M_{\pi}^{*}} 1 \bmod 2$ a well-defined invariant.
These counter-terms will depend on the normal operator $D^{\pi}_{A}$,
the Chern--Simons function and spectral invariants.

Note that in a ZHS the trivial orbit $\{[\mTheta]\}$ is always a
point in $\M^{r}_{\overline{\pi}}$ for every perturbation. In the
unperturbed case this is clear. In the presence of a perturbation
invariance by the stabilizer action at $\mTheta$ forces $\overline{\pi}_{\mTheta}=0$.

When $A=\mTheta$, the Dirac operator $D_{\mTheta}$ can be identified with the
canonical quaternionic linear Dirac operator on $S$ which we denote as $D$. The
operator $K_{A}$ (presently take $\pi=0$) is the boundary $B$
of the 4--dimensional signature operator, after identifying $\mLambda^{2}\cong\mLambda^{1}$
by the Hodge $*$--operator. To these  two operators $D$ and $B$ we
can associate the APS--spectral invariants \cite{APS}:
\begin{displaymath}
\eta(B),\quad \xi=\frac{1}{2}\Bigl( \eta(D)+\dimn_{\C}\kernel D\Bigr).
\end{displaymath}
If $X$ is compact oriented spin 4--manifold with oriented boundary $Y$
then an application of the APS index theorems to $X$ shows that
\begin{equation}\label{ind-equ}
\xi+\frac{1}{8}\eta(B) = -\mbox{Index}\, D^{(4)} -\frac{1}{8}\mbox{sign}\,X.
\end{equation}
Here $D^{(4)}$ is the Dirac operator on $X$ and ${\rm sign}\,X$ the signature.
Thus we see that the left-side of (\ref{ind-equ}) is always an integer.
As an aside, the $\bmod 2$ reduction of the right-side only involves the
signature term (since in four dimensions the Dirac operator is quaternionic linear and
so its index is even) and therefore is just the {\it Rokhlin invariant} $\mu(Y)$. 
Given a perturbation $\pi$ 
now set
\begin{displaymath}
c(g,\pi) =\xi+\frac{1}{8}\eta(B) +
\Bigl(\mbox{$\C$--spectral flow of 
$\{(1-t)D_{\mTheta}+t D^{\pi}_{\mTheta}\}_{t=0}^{1}$}
\Bigr)\in\Z.
\end{displaymath}
In the spectral-flow term $D_{\mTheta}$, $D^{\pi}_{\mTheta}$ are quaternionic linear
and thus $c(g,\pi)\equiv \mu(Y)\bmod 2$ continues to be true.
$c(g,\pi)$ is our counter-term associated to $\{[\mTheta]\}$.

\begin{remark}\rm Our convention for spectral flow is the 
the number of eigenvalues (counted algebraically) crossing $-\epsilon$ for $\epsilon>0$ sufficiently 
small. 
\end{remark}

In order to define the counter-terms associated with points in $\M^{r*}_{\overline{\pi}}$
we shall need two preliminaries. 
Firstly, consider the {\it normal spectral
flow} of $L^{\pi}_{A}$ along a path $\gamma$ in $\A$ i.e.
\begin{eqnarray*}
\SF^{\nu}(\gamma) &=& \mbox{\rm spectral flow of $D^{\pi}_{A}$ along
$\gamma$}
\end{eqnarray*}
which is defined because of Lemma~\ref{sf-fact-lem}.
On the reducible stratum $\A\subset\CC$, the Chern--Simons--Dirac
function reduces to the Chern--Simons function which we denote as $\cs$.
We remind the reader that $\cs$ depends on a basepoint which we choose to
be a trivial connection $\mTheta$ (which we fix once and for all).

\begin{lemma}\label{sf-calc}
Let $[x]$ be a point in $\B_{\A}$ and $[\gamma(t)]$, $t\in[0,1]$ a closed 
differentiable loop in $\B_{\A}$ based at $[x]$. Then
\begin{eqnarray*}
\SF^{\nu}(\gamma) &=& \cs(\gamma(1))-\cs(\gamma(0))\in \Z.
\end{eqnarray*}
\end{lemma}

\proof
First we invoke Lemma~\ref{complexification} which says we only need 
to compute the complex spectral flow for the complex Dirac operator
$D^{\C}_{A}$ on $S\otimes E=S\otimes_{\C}E$. According to \cite{APS}
this spectral flow coincides with the index of the four-dimensional Dirac
operator $D^{(4)}_{\widehat{A}}$ on the pull-back
$\widehat{S}\otimes\widehat{E}\to Y\times[0,1]$ of $S\otimes E\to Y$
with $\widehat{A}$ interpolating between $\gamma(0)$ at $Y\times\{0\}$
and $\gamma(1)$ at $Y\times\{1\}$. Since the initial and final 
connections are gauge equivalent, the boundary terms cancel in the 
application of the APS index theorem and we are left with
\begin{eqnarray*}
\phantom{Xi}\mbox{Index}\, D^{(4)}_{\widehat{A}} &=& -\int_{X}c_{2}(\widehat{A})\quad\mbox{\rm (Chern form)}
=-\frac{1}{8\pi^2}\int_{Y\times[0,1]}{\rm Tr}(F^2_{\widehat{A}})\\
&=& -\frac{1}{8\pi^2}\int_{Y\times[0,1]}
d\,{\rm Tr}\Bigl( a\wedge d_{\mTheta}a +\frac{2}{3}a \wedge a\wedge a\Bigr),
\quad a=\widehat{A}-\mTheta\\
&=& -\Bigl(\cs(\gamma(0))-\cs(\gamma(1))\Bigr).\hspace{2.1in}\qed
\end{eqnarray*}

The second preliminary: $\cs$ descends to a function $\overline{\cs}\colon 
\B_{\A}\to\R/\Z$ on the quotient space.
Since the value of $\overline{\cs}$ is constant on components of $\M^{r}$,
the image set $\overline{\cs}(\M^{r})$ is a finite number of values
$c_{1},\dots,c_{m}$ in $\R/\Z$. Let $\epsilon_{1}>0$ be the smallest
distance between pairwise distinct $c_{i}$'s where $\R/\Z$ has the distance
inherited from $\R$.
Let $\epsilon_{2}>0$ be the constant which is the smallest distance
between  pairwise distinct components of $\M^{r}$, in the metric
(\ref{metric-def}).

\begin{defn}\label{small-defn}\rm
Call a perturbation $\pi$ \textit{small} if 
$\overline{\cs}(\M^{r*}_{\overline{\pi}})$ is within an $\epsilon_{1}/3$--neighbourhood
of $\overline{\cs}(\M^{r*})$, and $\M^{r*}_{\overline{\pi}}$ is within an 
$\epsilon_{2}/3$--neighbourhood of $\M^{r*}$. 
\end{defn}

Assume $\pi$ to be 
small and non-degenerate in the sense of the preceding. This can be done by making
$\normi{\pi_{A,\mPhi}}{L^{2}_{2,A}}$
sufficiently small, by Proposition~\ref{adm-lem}. 
Write 
$$
\M^{r*}=\K^{1}\cup\dots\cup\K^{n}
$$ 
as the union of connected components. 
Then given any $[A]\in \M^{r*}_{\overline{\pi}}$ there is a
unique component $\K^{i}$ which is within $\epsilon_{1}/3$
of $[A]$. Denote by $\N^{i}$ the intersection of the
$\epsilon_{2}/3$--neighbourhood of $\K^{i}$ and the preimage under
$\overline{\cs}$ of the $\epsilon_{1}/3$--neighbourhood of $\overline{\cs}(\M^{r*})$ in
$\R/\Z$.
Let $[\gamma]$
be any path from $[\mTheta]$ to $[A]$. Let $[\overline{\gamma}]$
be any other path from $[\mTheta]$ to $\K^{i}$ with the property
that $[\gamma]$ and $[\overline{\gamma}]$ are homotopic 
relative to $\N^{i}\cup \{[\mTheta]\}$.
Then the expression
\begin{displaymath}
\kappa[A] = \SF^{\nu}(\gamma)
+\cs(\overline{\gamma}(0))-
\cs(\overline{\gamma}(1))\in\R
\end{displaymath}
is well-defined and independent of choice of $\gamma$ and $\overline{\gamma}$, 
by Lemma~\ref{sf-calc}.

Over $\B_{\A}$ we have the line bundle $\detind K^{\overline{\pi}}$ of the family 
of tangential operators $\{K^{\overline{\pi}}_{A}\}$. In contrast to
$\detind L^{\pi}$ this is an orientable line bundle. This is
basically the Taubes' orientation of $\M^{r*}_{\overline{\pi}}$
in \cite{taubes}. We fix the overall orientation by
specifying $\detind K^{\overline{\pi}}$ at $[\mTheta]$
by the following rule. The kernel and cokernel of $K^{\pi}_{\mTheta}$
are $\cong\su(2)$, the constant sections of $\ad E$, after $\ad E$ is
trivialized as $Y\times \su(2)$ by $\mTheta$. Orient 
$\detind K^{\overline{\pi}}_{\mTheta}$ as $o(\su (2))\wedge 
o(\su(2)^{*})$ where $o(\su(2))$ is a any chosen orientation and
$o(\su(2)^{*})$ is the dual orientation. We denote the induced
oriention at $[A]\in \M^{r*}_{\overline{\pi}}$ by 
$\epsilon[A]\in\{\pm 1\}$.

\begin{theorem}\label{thrm1}
Let {\it\,Y} be an oriented closed  integral
homology 3--sphere with  Riemannian
metric $g$. Let $\pi$ be a non-degenerate and small admissible perturbation for
$\M_{g,\pi}$, the perturbed Quaternionic Seiberg--Witten moduli
space with respect to $g$. 
The terms $c(g,\pi)$, $\epsilon[A]$ and $\kappa[A]$ as above are
well-defined and the sum
\begin{displaymath}
\tau(Y) = 
\sum_{[A]\in \M_{\pi}^{*}} 1 +\frac{1}{2}c(g,\pi) +
\sum_{[A]\in\M^{r*}_{\overline{\pi}}} \epsilon[A]\Bigl(\kappa[A]+\frac{1}{4}\Bigr)
\in \R\bmod 2\Z
\end{displaymath}
is independent of both $g$ and $\pi$ chosen. Furthermore $\tau(Y)$
does not depend on the orientation of {\it\,Y} and therefore $\tau$ defines an unoriented
diffeomorphism invariant for integral homology 3--spheres.
\end{theorem}

The extra term $1/4$ in the sum is inserted to make the invariant
independent of the orientation of $Y$.

Let $\lambda_{SU(3)}(Y)$ be the $SU(3)$--Casson invariant of Boden--Herald.
The definition of $\tau(Y)$ is modelled on $\lambda_{SU(3)}(Y)$
and both suffer from the defect that no multiple is obviously integral valued. This is
due to the usage of the Chern--Simons function. (Boden--Herald--Kirk \cite{BHK}
have devised an integer version of $SU(3)$--Casson 
that gets around the usage of Chern--Simons by an ad-hoc device.
It is not a completely natural
definition.)
However we have the following.

\begin{theorem}\label{thrm2}
Let {\it\,Y} be an integral homology 3--sphere and
$\lambda_{SU(3)}(Y)$ be the $SU(3)$--Casson invariant for $Y$. Then
\begin{displaymath}
\lambda_{SU(3)}(Y) + 2\tau(Y)
\end{displaymath}
is a $\Z\bmod 4\Z$--valued invariant of the unoriented diffeomorphism
type of $Y$.
\end{theorem}

The assertion of this theorem is that we have a cancellation
of the Chern--Simons terms, leaving only an integral expression.
Our contention is that combining $SU(3)$--Casson with an $SU(2)$--version
of Seiberg--Witten is the natural way of presenting the topological
information contained in the two theories. This will be worked out in 
greater detail in a further article where a unified 
approach to the two theories and an
integer valued Seiberg--Witten/Casson invariant is defined.

The proof of the Theorem~\ref{thrm1} is in \S\ref{proof1-sect} and
Theorem~\ref{thrm2} in \S\ref{proof2-sect}.

\section{Compactness}\label{comp-sect}

In this section we prove Proposition~\ref{adm-lem}.
Recall that our 3--manifold $Y$ is assumed to be Riemannian with metric $g$.
We shall need to vary $g$ at two points in this article.  
In the present section we shall utilize rescaling $g$ to
establish compactness of the moduli space. In \S\ref{proof1-sect}
we shall analyse the change in the moduli space as $g$ varies
in a 1--parameter family.

We set-up a framework for comparing the SW--equation
for different metrics. Spinors and in particular the Dirac operator
are not canonically associated objects to a Riemannian structure.

The first task is to fix a model for the spin structure and spinors.
Our metric $g$ shall be taken as the reference.
On a compact
3--manifold we can always find a smooth nowhere vanishing
vector field, let us denote this as $e_1$. Additionally assume it
is of unit length with respect to $g$. By working perpendicular
to $e_1$ we can complete this to a global orthonormal
frame $(e_1,e_2,e_3)$. Assume the orientation $e_1\wedge e_2\wedge e_3$
coincides with the orientation on $Y$.
This global frame defines a trivialization $Y\times SO(3)$ of the (positively)
oriented orthonormal frame bundle of $Y$.

Let $sp(1)={\rm spin}(3)\subset CL(Y)$ denote the unit quaternions and fix a group homomorphism
$sp(1)\to SO(3)$ which is the 2--fold covering map. Then we fix the
spin structure on $Y$ (with respect to $g$) by the projection
\begin{displaymath}
P=Y\times sp(1) \to Y\times SO(3).
\end{displaymath}
The spinor bundle $S$ is then given by $P\times_{\varrho}\qu$ where $\varrho$ is the fundamental
representation of $sp(1)$ on $\qu$. Since $[(y,h),q] = [(y,1),hq]$, $S$  has a natural
trivialization as $Y\times\qu$ and sections of $S$ are simply the $\qu$--valued 
functions on $Y$. Notice that
the quaternionic structure on $S$ is exactly right multiplication
on the $\qu$ factor of $Y\times\qu$.

The trivialization $Y\times SO(3)$ of $TY$ also induces a trivialization
$Y\times CL(\R^3)$ of the Clifford bundle $CL(Y)$ with the constant section
$\widehat{e}_1=(1,0,0)$ corresponding to the vector field $e_1$, $\widehat{e}_2=(0,1,0)$ to
$e_2$ etc.
Fix the (left) Clifford representation on $\qu$ of the Clifford
algebra $CL(\R^3)$ by mapping
\begin{displaymath}
\widehat{e}_1\mapsto i,\quad \widehat{e}_2\mapsto j,
\quad \widehat{e}_3\mapsto k.
\end{displaymath}
That is to say, $\widehat{e}_{i}\cdot h = ih$ etc.
On $S\otimes_{\qu}E$ the Dirac operator
now takes the form
\begin{displaymath}
D_{A} =(i\otimes 1)\nabla^{A}_{e_{1}} +(j\otimes 1)\nabla^{A}_{e_{2}} 
+(k\otimes 1)\nabla^{A}_{e_{3}}.
\end{displaymath}
Suppose now we want to change the metric from $g$. This is
achieved by pulling back the metric $g$ by an automorphism
$h$ of $TY$. Using the frame $(e_1,e_2,e_3)$ as a basis
can conveniently think of $h$ as a smooth map $h\colon Y\to GL(3)$.
The global frame $(e_1,e_2,e_3)$
is pulled back to a global frame $(h^{-1}(e_1),h^{-1}(e_2),h^{-1}(e_3))$
for the pulled back metric. In the same way as above this global frame
defines a trivialization
$Y\times SO(3)$ of the oriented orthonormal frame bundle in the pulled 
back metric and we may proceed with the spin structure, spinors 
etc.\ as constructed before. In particular we notice that the model for the
spinor bundle as $\qu$--valued functions on $Y$ remains the same
in the pulled back metric but the Clifford mutiplication
changes and is now defined by
\begin{displaymath}
\widehat{h^{-1}(e_1)}\mapsto i,\quad
\widehat{h^{-1}(e_2)}\mapsto j, \quad
\widehat{h^{-1}(e_3)}\mapsto k.
\end{displaymath}
If $h$ is actually an isometry with respect to
$g$ then we are merely changing the trivialization of $S$.

Let $g'$ denote the new metric defined by $h$ and $\nabla^{g'}$ the spin
connection on $S$. Then the Dirac operator coupled to $A$ with
respect to $g'$ is given by
\begin{displaymath}
D^{g'}_A = (i\otimes 1)\nabla^{g',A}_{h^{-1}(e_1)} + (j\otimes 1)\nabla^{g',A}_{h^{-1}(e_2)} 
+(k\otimes 1)\nabla^{g',A}_{h^{-1}(e_3)}.
\end{displaymath}
Similiarly one may obtain expressions for the bilinear forms 
$\{\cdot\}_{0}$ and $B$ with respect to
$g'$ in terms of $h$.

Consider now the special case when $g$ is rescaled as 
$g_{\lambda}=\lambda^2 g$ where $\lambda>0$ is a constant.
Clearly $g_{\lambda}$ is induced by $h= \lambda {\rm Id} $ so 
$h^{-1}(e_i) = e_i/\lambda$. Under the above model for the spinors, the Hermitian
metric on $S$ is fixed. However, we may choose to vary this with $\lambda$.
In the present case, for $g_{\lambda}$ we may set 
\begin{equation}\label{spinor-met}
\langle{\cdot,\cdot}\rangle_{\lambda}=\lambda^{\alpha}\langle{\cdot,\cdot}\rangle
\end{equation}
where the right-hand inner product is the original one on $S$. A good
choice for $\alpha$ will be made later. In the next lemma, a \lq$\lambda$\rq\ superscript
means an object taken with respect to the metric $g_{\lambda}$. Unmarked
objects are taken with respect to $g$.

\begin{lemma}
Fix the model for spinor bundle $S$ by $g$, and use the spinor metric given by
(\ref{spinor-met}) in the Riemannian metric $g_{\lambda}$. Then the
following hold.
\begin{enumerate}
\item $D^{\lambda}_A = (1/{\lambda}) D_A$ \label{eq-d-1}
\item $\{\cdot\}_0^{\lambda} = \lambda^{2+\alpha} \{\cdot\}_0$ \label{eq-d-2}
\item $d^{*,\lambda}_A = (1/{\lambda^2})d^*_A$ on $\mOmega^1\otimes \ad E$\label{eq-d-3}
\item $B^{\lambda} = \lambda^{\alpha} B$\label{eq-d-4}
\item $*^{\lambda} = (1/{\lambda})*$ on $\mLambda^2$\label{eq-d-5}
\end{enumerate}
\end{lemma}

\proof For (\ref{eq-d-1}) recall that the Levi--Civita connection is invariant
under rescaling the metric by a constant. This leaves the connection term 
$\nabla^{g',A}=\nabla^{g,A}$. The formula now follows from $h^{-1}(e_{i})=e_{i}/\lambda$.
For (\ref{eq-d-2}) establish the rule 
$\omega\cdot_{\lambda}\phi=\frac{1}{\lambda^{2}}\omega\cdot\phi$ and
$\alpha\cdot_{\lambda}\phi=\frac{1}{\lambda}\alpha\cdot\phi$ where $\omega$
is a 2--form and $\alpha$ a 1--form. The new coframe $e^{\lambda *}_{i}=\lambda e^{*}_{i}$
and so the action of $e^{*}_{i}$ with respect to $g_{\lambda}$ is $1/\lambda$ of
the action with respect to $g$. For (\ref{eq-d-3}) in the defining
equation $\int\langle{d_{A}\gamma,a}\rangle_{\lambda}dg_{\lambda}=
\int\langle{\gamma,d^{*,\lambda}_{A}a}\rangle_{\lambda}dg_{\lambda}$,
we have
$\langle{\gamma,d^{*,\lambda}_{A}a}\rangle_{\lambda}=
\lambda^{-2}\langle{\gamma,d^{*,\lambda}_{A}a}\rangle$.
For (\ref{eq-d-4}) the defining equation is $\langle{\gamma(\phi),\psi}\rangle_{\lambda}
=\lambda^{\alpha}\langle{\gamma(\phi),\psi}\rangle=
\lambda^{\alpha}\langle{\gamma, B(\phi,\psi)}\rangle=
\langle{\gamma,B^{\lambda}(\psi,\phi)}\rangle$. (\ref{eq-d-5}):
$*(e_{1}^{*}\wedge e_{2}^{*})=e_{3}^{*}$ and $*^{\lambda}(\lambda e_{1}^{*}\wedge
\lambda e_{2}^{*}) =\lambda e_{3}^{*}$, etc. \endproof

The preceding lemma easily implies the following principle result we need
on rescaling the metric:

\begin{prop}\label{rescale-prop}
Fix the model for the spinor bundle $S$ by $g$, and let $S$ have the
metric (\ref{spinor-met}) with respect to $g_{\lambda}$ where $\alpha=-2$. Then 
the perturbed SW--equation (\ref{sw-eq2}) with respect to $g$ is equivalent to the 
following equation with respect to $g_{\lambda}$:
\begin{equation} \label{sw-eq3}
\left\{
\begin{array}{ccl}
\displaystyle F_A - \{\mPhi\cdot\mPhi\}_0^{\lambda}
+k_{A,\mPhi}&=&0\\
\noalign{\vskip 5pt}
\displaystyle D^{\lambda}_A {\mPhi} 
+\frac{1}{\lambda}l_{A,\mPhi} &=& 0.
\end{array}
\right.
\end{equation}
Furthermore the perturbation $\pi^{\lambda}(A,\mPhi)=
(*^{\lambda}k_{A,\mPhi},
(1/\lambda)l_{A,\mPhi} )$ is an admissible
perturbation with respect to $g_{\lambda}$.
\end{prop}

The scheme of the proof of the compactness of the moduli space
rests on a Bochner argument to get a $L^4$--bound
on the spinors, Uhlenbeck's Theorem \cite{uhlen}
and as mentioned above, rescaling. In the 4--dimensional context such
an argument is presented in Feehan--Leness \cite{FL}.
The basic input is contained in the following two lemmas.

\begin{lemma}\label{bound-lem}
Let $(A,\mPhi)$ be a solution of the perturbed SW--equation (\ref{sw-eq3}),
defined on {\it\,Y} with respect to the metric $g_{\lambda}=\lambda^2g$. 
Let $s$ denote the scalar curvature of {\it\,Y} with respect to $g$.
Then
\begin{eqnarray*}
\int_{Y}|\mPhi|^{4}_{\lambda}dg_{\lambda}
&\le& 
\frac{8}{\lambda}\int_{Y}\frac{s^{2}}{16}+|k_{A,\mPhi}|^{2}+|l_{A,\mPhi}|^{2} dg\\
\int_{Y}|k_{A,\mPhi}|^{2}_{\lambda}dg_{\lambda} 
&\le&
\frac{1}{\lambda}\int_{Y}|k_{A,\mPhi}|^{2}dg\\
\int_{Y}\Bigl|\frac{1}{\lambda} l_{A,\mPhi}\Bigr|^{2}_{\lambda}dg_{\lambda}
&\le&
\frac{1}{\lambda}\int_{Y}|l_{A,\mPhi}|^{2}dg.
\end{eqnarray*}
The spinor metric (\ref{spinor-met}) on the left-side is taken with $\alpha=-2$.
\end{lemma}

\proof This is a straightforward manipulation
involving the Bochner formula for the Dirac operator which reads:
\begin{displaymath}
(D^{\lambda}_A)^* D^{\lambda}_A\mPhi = (\nabla^{\lambda}_A)^* \nabla^{\lambda}_A\mPhi 
+\frac{1}{4} s^{\lambda}\mPhi +F_A\cdot_{\lambda}\mPhi.
\end{displaymath}
Here and below a \lq$\lambda$\rq\ subscript or superscript 
indicates the object taken with respect to $g_{\lambda}$.
Unscripted objects are taken with respect to $g$.
Taking the inner product with $\mPhi$ and integrating gives
$$
\int_Y |D^{\lambda}_A\mPhi|^2_{\lambda}\,dg_{\lambda}
=\int_Y |\nabla^{\lambda}_A\mPhi|^2_{\lambda}\,dg_{\lambda}
+\int_Y \frac{1}{4}s^{\lambda}|\mPhi|^2_{\lambda}\,dg_{\lambda}
+\int_Y \langle{F_A,\{\mPhi\cdot \mPhi\}^{\lambda}_{0}}\rangle_{\lambda}\,dg_{\lambda}.
$$
Applying the SW--equation (\ref{sw-eq3}) and after some manipulation we obtain
\begin{displaymath}
\int_Y \biggl( \frac{s^{\lambda}}{4}
 -|k_{A,\mPhi}|_{\lambda}\biggr)|\mPhi|^2_{\lambda} \,dg_{\lambda}
+\frac{1}{2}\int_Y |\mPhi|^4_{\lambda}\,dg_{\lambda}
\leq \frac{2}{\lambda^2} \int_Y |l_{A,\mPhi}|^2_{\lambda}\,dg_{\lambda}.
\end{displaymath}
This in turn implies
\begin{displaymath}
\int_Y |\mPhi|^4_{\lambda}\,dg_{\lambda} \leq  
2\Gamma_{\lambda}\,\biggl(\, \int_Y |\mPhi|^4_{\lambda}\,dg_{\lambda}\biggr)^{1/2}
+\frac{4}{\lambda^2}\int_Y |l_{A,\mPhi}|^2_{\lambda}\,dg_{\lambda}
\end{displaymath}
where $\Gamma_{\lambda}\ge 0$ is given by
$$
\Gamma^2_{\lambda} = \int_Y\biggl(\,\frac{s^{\lambda}}{4}
-|k_{A,\mPhi}|_{\lambda}\biggr)^2dg_{\lambda}.
$$
Therefore
\begin{equation}\label{bd-eq1}
\int_Y |\mPhi|^4_{\lambda}\,dg_{\lambda} 
\leq  4\Gamma^2_{\lambda} + 
\frac{8}{\lambda^2}\int_Y |l_{A,\mPhi}|^2_{\lambda}\,dg_{\lambda}.
\end{equation}
Under rescaling the metric from $g$ to $g_{\lambda}=\lambda^2g$ we have
$dg_{\lambda}=\lambda^{3} dg$ and the following relations
hold:
\begin{eqnarray}
\int_Y (s^{\lambda})^2\, dg_{\lambda} 
&=& \int_Y (\lambda^{-2}s)^2\lambda^3 \,dg \nonumber\\
\int_Y|k_{A,\mPhi}|^2_{\lambda}\,dg_{\lambda} 
&=& \int_Y\lambda^{-4}|k_{A,\mPhi}|^2 \lambda^3\,dg\label{bd-scale}\\
\int_Y |l_{A,\mPhi}|^2_{\lambda}\,dg_{\lambda}
&=& \int_Y \lambda^{-2}|l_{A,\mPhi}|^2\lambda^3\,dg.\nonumber
\end{eqnarray}
Hence
\begin{displaymath}
\Gamma^2_{\lambda} 
\leq \frac{4}{\lambda}\int_Y 
\Bigl(\frac{s^2}{16} +|k_{A,\mPhi}|^2\Bigr)dg
\end{displaymath}
Together with
(\ref{bd-eq1}) and (\ref{bd-scale}) we get the desired bounds. \endproof

Introduce the notation $B_r$ for the closed Euclidean ball of radius $r$ in $\R^3$.
Fix a model for the spinors $S$ on $B_1$ with respect to the Euclidean metric
as in the preceding and let $E_{0}=B_1\times \C^{2}$ denote the trivial
$SU(2)$--bundle. This trivialization defines the canonical trivial connection
$d$ on $E$. 

\begin{lemma}\label{basic-lem}
Allow any metric on $B_{1}$.
Let the pair $(A=d+a,\mPhi)\in L^{2}_{2}$ be defined on $E_{0}\to B_{1}$.
Assume that (a) $d^{*}a=0$ (b) $\normi{a}{L^{2}_{1}}\le C_{1}$, 
$\normi{\mPhi}{L^{4}}\le C_{2}$ 
(c) $(A,\mPhi)$ satisfies a perturbed SW--equation of the form 
(\ref{sw-eq2}) on $E$ with $\pi_{A,\mPhi}=(*k,l)_{A,\mPhi}\in L^{2}_{2}$, and
(d) $\normi{\pi_{A,\mPhi}}{L^{2}_{2}}\le C_{3}$. Then
$\normi{a}{L^2_3(B_{1/2})}$, $\normi{\mPhi}{L^2_3(B_{1/2})}$ are uniformly bounded
independent of $a$ and $\mPhi$.
\end{lemma}

\proof
We may rewrite the equations both $a$ and $\mPhi$ satisfy as
\begin{eqnarray*}
(d+d^{*})a &=& -a\wedge a +\{\mPhi\cdot\mPhi\}_{0} -k_{A,\mPhi}\\
D\mPhi &=& -a\cdot\mPhi -l_{A,\mPhi}.
\end{eqnarray*}
Here $D$ is the canonical Dirac operator associated with $B_{1}$ tensored with the trivial
factor $\C^{2}$. 
$\normi{a}{L^{4}}$ is uniformly bounded by the Sobolev embedding
$L^{2}_{1}\subset L^{4}$ and condition (a). 
The terms $k_{A,\mPhi}$, $l_{A,\mPhi}$ being uniformly bounded in 
$L^{2}_{2}$ are uniformly bounded in $C^{0}$. Since 
$\normi{a\cdot\mPhi}{L^{2}}\leq \normi{a}{L^{4}}\normi{\mPhi}{L^{4}}$ we see that
$D\mPhi$ is uniformly bounded in $L^{2}$. The
basic elliptic inequality 
$\normi{\mPhi}{L^p_{k+1}(B_{r'})}\leq\const(\normi{D\mPhi}{L^{p}_{k}(B_{r})}+
\normi{\mPhi}{L^{p}(B_{r})})$, $r'\le r$
for $D$ forces $\mPhi$ to be uniformly bounded in $L^{2}_{1}$ over
$B_{r_{1}}$, $r_{1}<1$. The embedding $L^{2}_{1}\subset L^{6}$ now makes both $a$ and
$\mPhi$ uniformly bounded in $L^{6}$ over $B_{r_{1}}$. The bound
$\normi{a\cdot\mPhi}{L^{3}(B_{r_{1}})}\leq \normi{a}{L^{6}(B_{r_{1}})}
\normi{\mPhi}{L^{6}(B_{r_{1}})}$ now makes $D\mPhi$ uniformly bounded in $L^{3}(B_{r_{1}})$ and
thus $\mPhi$ is uniformly bounded in $L^{3}_{1}(B_{r_{2}})$, $r_{2}<r_{1}$ and therefore
$L^{p}(B_{r_{2}})$, $2\le p<\infty$. Now observe 
$\normi{a\cdot\mPhi}{L^{4}(B_{r_{2}})}\leq \normi{a}{L^{6}(B_{r_{2}})}
\normi{\mPhi}{L^{12}(B_{r_{2}})}$ and by repeating the argument we get $\mPhi$
uniformly bounded in $L^{4}_{1}(B_{r_{3}})$, $r_{3}<r_{2}$. A similiar type of
argument using the elliptic estimate for $d+d^{*}$ also 
establishes that $a$ is uniformly bounded in $L^{4}_{1}(B_{r_{3}})$.

To obtain uniform bounds for $a$ and $\mPhi$ in $L^{2}_{2}(B_{r_{4}})$, $r_{4}<r_{3}$
we need to obtain uniform bounds for the quadratic terms $a\wedge a$, $\{\mPhi\cdot\mPhi\}_{0}$
and $a\cdot\mPhi$ in $L^{2}_{1}(B_{r_{3}})$. However this follows from the continuous
multiplication $L^{4}_{1}(B_{r_{3}})\times L^{4}_{1}(B_{r_{3}})\to L^{2}_{1}(B_{r_{3}})$. Finally
this puts $a$ and $\mPhi$ in the continuous range for Sobolev multiplication and from this
a uniform bound in $L^{2}_{3}(B_{r_{5}})$, $r_{5}<r_{4}$ is obtained. \endproof

\begin{prop}\label{compact}
$\M_{\pi}$ is a compact subspace of $\B$ where $\pi$ an admissible perturbation.
That is to say, given any sequence $(A_i,\mPhi_i)$ of $L^2_2$--solutions
to (\ref{sw-eq2}) there is a subsequence $\{ i'\}\subset\{i\}$
and $L^2_3$--gauge transformations $g_{i'}$ such that $g_{i'}(A_{i'},\mPhi_{i'})$
converges in $L^2_2$ to a solution of the $\pi$--perturbed SW--equations.
\end{prop}

\proof
By Proposition~\ref{rescale-prop}
a solution $(A,\mPhi)$ of (\ref{sw-eq2})
is equivalent to a solution of (\ref{sw-eq3}), the SW--equation with respect to
$g_{\lambda}$ and with perturbation $\pi^{\lambda}$.
Thus it suffices to prove compactness
of the moduli space $\M_{g_{\lambda},\pi^{\lambda}}$ of solutions of (\ref{sw-eq3}) for
any $\lambda>0$.

Choose $\lambda$ large such that any 
geodesic ball $B$ of unit radius in $Y$ is sufficiently close to the Euclidean
metric in $C^{3}$, so that Uhlenbeck's Theorem \cite{uhlen} applies over $B$. 
Let $\epsilon_{0}>0$ be the constant in Uhlenbeck's Theorem such that if
any $L^{2}_{1}$ connection $A$ on $E|_{B}$ satisfies $\normi{F_{A}}{L^{2}(B)}<\epsilon_{0}$
then there is a gauge transformation $g\in L^{2}_{2}(B)$ which changes $A$ so that 
$g(A)=d+a$ is in Coloumb gauge  $d^{*}a=0$ and $\normi{a}{L^{2}_{1}(B)}\le c 
\normi{F_{A}}{L^{2}(B)}$. Here we use a fixed trivialization
$E|_{B}\cong B\times\C^{2}$ with trivial connection $\nabla$ or $d$.

Assume that $(A,\mPhi)$ is a solution of (\ref{sw-eq3}).
The proof of Lemma~\ref{basic-lem} gives us an additional fact. It
shows that $a$ is of class $L^{2}_{2}(B_{1/2})$ and by a straightforward
bootstrapping argument we see that $g$ is actually in $L^{2}_{3}(B_{1/2})$. 

In the definition of an admissible perturbation
$\normi{\pi^{\lambda}_{A,\mPhi}}{L^{2}_{2,A}}$ is uniformly bounded for every $\lambda>0$.
In order to apply Lemma~\ref{basic-lem} we need to deduce a uniformly bound for
$\normi{\pi^{\lambda}_{A,\mPhi}}{L^{2}_{2}(B)}$. The covariant derivatives $\nabla$
and $\nabla^{A}$ upto second order are related by
\begin{eqnarray*}
\nabla\omega &=& \nabla^{A}\omega - a(\omega) \\
\nabla^{2}\omega 
&=& (\nabla^{A}-a)(\nabla^{A}\omega -a (\omega))\\
&=& (\nabla^{A})^{2}\omega +(\nabla a)(\omega) +2a(a(\omega)).
\end{eqnarray*}
Utilizing the embedding $L^{2}_{2}(B)\subset C^{0}(B)$ and $L^{2}_{1}(B)\subset L^{4}(B)$ we obtain
\begin{eqnarray*}
\normi{\nabla\omega}{L^{2}(B)} 
&\le & \const\Bigl(\normi{\nabla^{A}\omega}{L^{2}(B)} 
+\normi{a}{L^{2}(B)}\normi{\omega}{L^{2}_{2}(B)}\Bigr)\\
\normi{\nabla^{2}\omega}{L^{2}(B)}
&\le & \const\Bigl(\normi{(\nabla^{A})^{2}\omega}{L^{2}(B)} 
+ \normi{a}{L^{2}_{1}(B)}\normi{\omega}{L^{2}_{2}(B)}\\
&&\hspace{2.2in}+\normii{a}{L^{2}_{1}(B)}{2}\normi{\omega}{L^{2}_{2}(B)}\Bigr).
\end{eqnarray*}
Choose $\epsilon_{1}\le\epsilon_{0}$ so that 
$\normi{F_{A}}{L^{2}(B)}<\epsilon_{1}$
forces $\normi{a}{L^{2}_{1}(B)}$ to be very small; then
the error terms $\normi{a}{L^{2}(B)}\normi{\omega}{L^{2}_{2}(B)}$, 
$\normi{a}{L^{2}_{1}(B)}\normi{\omega}{L^{2}_{2}(B)}$  and 
$\normii{a}{L^{2}_{1}(B)}{2}\normi{\omega}{L^{2}_{2}(B)}$ are $\ll \normi{\omega}{L^{2}_{2}(B)}$
and we get a uniform estimate 
$\normi{\omega}{L^{2}_{2}(B)}\le \const\normi{\omega}{L^{2}_{2,A}(B)}$.

Lemma~\ref{bound-lem} shows that $\normi{\mPhi}{L^{4}}$ is uniformly bounded with 
respect to $g_{\lambda}$ and $\normi{F_{A}}{L^{2}}\to 0$ as $\lambda\to\infty$.
Increase $\lambda$ if necessary so that 
$\normi{F_{A}}{L^{2}(B)}<\epsilon_{1}$ for all $B$. 
Suppose now that $(A_{i},\mPhi_{i})$ is a sequence of solutions of (\ref{sw-eq3}). 
Denote by $B_{1/2}$ the geodesic ball with the same center as $B$ but half the radius. 
Uhlenbeck's Theorem and the uniform bounds of
Lemma~\ref{basic-lem} finds $L^{2}_{3}$ gauge transformations $g_{i}$ over $B_{1/2}$
such that after passing to a subsequence, $g_{i}(A_{i},\mPhi_{i})$
converges in $L^{2}_{2}(B_{1/2})$ to a SW--solution (\ref{sw-eq3}) over $B$.
Now the standard covering argument in \cite[\S 4.4.2]{DK} (also see \cite{FL}) shows
that after global gauge transformations and passing to subsequences, $(A_{i},\mPhi_{i})$
can be made to converge in $L^{2}_{2}$ over all of $Y$. \endproof

The preceding proof also shows:

\begin{cor}\label{regular}
Let $(A,\mPhi)$ be a perturbed SW--solution (\ref{sw-eq2}). There is an $L^{2}_{3}$
gauge transformation $g$ such that $g(A,\mPhi)$ is in $L^{2}_{3}$.
\end{cor}

\begin{cor}
There is an $\epsilon_{0}>0$ such that for any 
$0<\epsilon<\epsilon_{0}$,
if $\normi{\pi_{A,\mPhi}}{L^{2}_{2,A}}\!<\epsilon$
uniformly then given any $[A,\mPhi]\in \M_{\pi}$ there is a $[A',\mPhi']\in\M$ such
that $d([A,\mPhi], [A',\mPhi'])<\epsilon$, $d$ being the metric (\ref{metric-def}).
\end{cor}

\proof Suppose false. Then there exists sequences $\{\pi_{i}\}$ and
$\{(A_{i},\mPhi_{i})\}$ with $[A_{i},\mPhi_{i}]\in \M_{\pi_{i}}$  such that 
$\normi{(\pi_{i})_{A_{i},\mPhi_{i}}}{L^{2}_{2,A_{i}}}\to 0$
but with $d([A_{i},\mPhi_{i}],[A',\mPhi'])$ bounded away from zero over
$[A',\mPhi']\in\M$.
The sequence also satisfies
\begin{equation}\label{lim-equ}
\normi{F_{A_{i}}-\{\mPhi_{i}\cdot\mPhi_{i}\}_{0}}{L^{2}}
+\normi{D_{A_{i}}\mPhi_{i}}{L^{2}}\to 0.
\end{equation}
The proof of Proposition~\ref{compact} shows that after gauge transformations and passing
to a subsequence which we shall also denote as $(A_{i},\mPhi_{i})$, 
$(A_{i},\mPhi_{i})$ converges in $L^{2}_{2}$ and the limit, by
(\ref{lim-equ}) is necessarily a unperturbed SW--solution. This
is a contradiction. \endproof

\section{Construction of perturbations}\label{pert-sect}

In this section we prove Proposition~\ref{adm2-lem}.
Introduce the notation
$B(\epsilon)$ for the $\epsilon$--ball in the slice space
$X_{A,\mPhi}$. (Recall this is a Hilbert space in an $L^{2}_{2}$--Sobolev norm.)
Denote by $\beta\colon X_{A,\mPhi}\to[0,1]$ a smooth cut-off function
with support in $B(\epsilon)$.

\begin{lemma}\label{implicit}
Fix $(A,\mPhi)$. There is an $\epsilon>0$ and a differentiable function
$\xi\colon B(\epsilon)\times X_{A,\mPhi}\to (\kernel \delta^{0}_{A,\mPhi})^{\perp}
\subset L^{2}_{3}( \ad E)$ such that given any 
$(a,\phi;b,\psi)\in B(\epsilon)\times X_{A,\mPhi}$, the relation
\begin{equation}\label{xi-equ}
(b,\psi) +\delta^{0}_{A,\mPhi}\xi (a,\phi;b,\psi) \in X_{A+a,\mPhi+\phi}
\end{equation}
holds. Here $(\kernel \delta^{0}_{A,\mPhi})^{\perp}$ denotes the 
$L^{2}$--orthogonal complement.
\end{lemma}

\proof Apply the Implicit Function theorem 
to the map
\begin{displaymath}
H(\xi,(a,\phi),(b,\psi)) =\delta^{0*}_{A+a,\mPhi+\phi}\delta^{0}_{A,\mPhi}(\xi)
+\delta^{0*}_{A+a,\mPhi+\phi}(b,\psi)
\end{displaymath}
from $(\kernel \delta^{0}_{A,\mPhi})^{\perp}\times B(\epsilon)\times X_{A,\mPhi}
\to (\kernel \delta^{0}_{A,\mPhi})^{\perp}\cap L^{2}_{1}$. 
The linearization of $H$ at $(0,0,0)$ restricted to 
$(\kernel \delta^{0}_{A,\mPhi})^{\perp}$ is an isomorphism. This establishes
the existence of the function $\xi=\xi(a,\phi;b,\psi)$ but only for $(a,\phi)$ and $(b,\psi)$
defined in sufficiently small neighbourhoods of zero. However notice that
if $(b,\psi)$ satisfies (\ref{xi-equ}) then 
for any real constant $c$,
$c(b,\psi)$ satisfies the same equation
but with $\xi$ replaced by $c\xi$. That is we can allow the $(b,\psi)$
to be defined in $\xi$ for all $X_{A,\mPhi}$ by extending $\xi$ linearly in that factor.
\endproof

Let us now assume $\mPhi\neq 0$. Set $\epsilon>0$ to  be less than the 
constant  in Lemma~\ref{implicit} and also such that
$B(\epsilon)$ injects into $\B$. Assume $\supp\beta\subset B(\epsilon)$.
Fix $(b,\psi)\in X_{A,\mPhi}$.
Define a function $(A,\mPhi)+X_{A,\mPhi}\to L^{2}_{2}(\mLambda^{1}\otimes \ad E)\times
L^{2}_{2}(S\otimes_{\qu}E)$ by the rule
\begin{equation}\label{pi-equ}
\pi_{A+a,\mPhi+\phi} = \beta(a,\phi)(b,\psi)
+\delta^{0}_{A,\mPhi}\xi(a,\phi;\beta(a,\phi)(b,\psi))
\end{equation}
for $(a,\phi)\in B(\epsilon)$. By construction $\pi$ has
support in $B(\epsilon)$.
Extend $\pi$ to
$\CC$ by $\G$--equivariance.
Clearly $\pi_{A+a,\mPhi+\phi}\in X_{A+a,\mPhi+\phi}$ and 
$\pi_{A,\mPhi}=(b,\psi)$.

\begin{lemma}\label{pert1-lem}
For $\epsilon>0$ sufficiently small, the perturbation $\pi$ in (\ref{pi-equ})
satisfies a uniform bound
$\normi{\pi_{A,\mPhi}}{L^{2}_{2,A}} \le C$.
\end{lemma}

\proof
$\xi$ satisfies $H(\xi,(a,\phi),\beta.(b,\psi))=0$.
Thus
\begin{equation}\label{D-equ}
\Delta_{A,\mPhi}\xi +N_{1}(a,\phi)\xi + N_{2}(a,\phi)(b,\psi) + 
\delta^{0*}_{A,\mPhi}(\beta.(b,\psi))=0
\end{equation}
where $\Delta_{A,\mPhi}$ is
a second order elliptic operator with coefficients depending on $(A,\mPhi)$
and $N_{1}$ and $N_{2}$ are lower order terms. $N_{1}$ is a bilinear expression
in $(a,\phi)$ and $\delta^{0}_{A,\mPhi}(\xi)$. $N_{2}$ is a bilinear expression
in $(a,\phi)$ and $(b,\psi)$.
After some calculation it is seen that
$N_{1}$, $N_{2}$ satisfy, by Sobolev theorems
\begin{eqnarray}
\normi{N_{1}(a,\phi)\xi}{L^{2}_{1}}&\leq& \const \normi{(a,\phi)}{L^{2}_{2}}
\normi{\xi}{L^{2}_{3}}\label{N-equ}\\
\normi{N_{2}(a,\phi)(b,\psi)}{L^{2}_{1}} &\leq &\const
\normi{(a,\phi)}{L^{2}_{2}}\normi{(b,\psi)}{L^{2}_{2}}.\nonumber
\end{eqnarray}
On the other hand since $\Delta_{A,\mPhi}$ is invertible 
on $(\kernel\delta^{0}_{A,\mPhi})^{\perp}$,
\begin{equation}\label{D2-equ}
\normi{\xi}{L^{2}_{3}} \leq\const\normi{\Delta_{A,\mPhi}\xi}{L^{2}_{1}}.
\end{equation}
Now make $\epsilon>0$ sufficiently
small so that $\normi{(a,\phi)}{L^{2}_{2}}$ is correspondingly small.
Then
(\ref{D-equ}), (\ref{N-equ}) and (\ref{D2-equ}) give
$\normi{\xi}{L^{2}_{3}} \leq\const\normi{(b,\psi)}{L^{2}_{2}}$.
Thus by (\ref{pi-equ}) we have a uniform bound
\begin{displaymath}
\normi{\pi_{A+a,\mPhi+\phi}}{L^{2}_{2}}
\leq \const\normi{(b,\psi)}{L^{2}_{2}}\leq C.
\end{displaymath}
In the above the Sobolev norms were taken with respect to some fixed connection
$A_{0}$, which is commensurate to the Sobolev norm taken to say $A$.
If $\normi{a}{L^{2}_{2}}$ is sufficiently small then 
\begin{eqnarray*}
\normi{\nabla^{A+a}\pi_{A+a,\mPhi+\phi}}{L^{2}}&\leq& \const
\normi{\nabla^{A}\pi_{A+a,\mPhi+\phi}}{L^{2}},\\
\normi{\nabla^{A+a}\nabla^{A+a}\pi_{A+a,\mPhi+\phi}}{L^{2}}&\leq& \const
\normi{\nabla^{A}\nabla^{A}\pi_{A+a,\mPhi+\phi}}{L^{2}}
\end{eqnarray*}
uniformly. By reducing
$\epsilon$ again if necessary,
the bound 
$\normi{\pi_{A,\mPhi}}{L^{2}_{2,A}} \le C$ is established. \endproof

This lemma directly shows

\begin{prop}\label{pert-1}
Assume $\mPhi\neq 0$.
Given any $(b,\psi)\in X_{A,\mPhi}$ there is an admissible perturbation $\pi$
such that $\pi_{A,\mPhi}=(b,\psi)$. Furthermore the support of $\pi$
may be chosen to be contained in an arbitarily small $\G$--invariant 
neighbourhood of the orbit  $\G\cdot(A,\mPhi)$.
\end{prop}

The slice at a reducible $(A,0)$
has a natural splitting $X_{A,0}=X^{r}_{A}\times L^{2}_{2}(S\otimes_{\qu}E)$.
The stabilizer of $(A,0)$ (which is $\{\pm 1\}$, $U(1)$ or $SU(2)$) acts diagonally on both 
of the factors
$X^{r}_{A}$ and $L^{2}_{2}(S\otimes_{\qu}E)$. If $\pi$ is a perturbation
then the stabilizer action forces the normal or spinor component of $\pi_{A,0}$ to
be zero, since $\pi$ is required to be $\G$--equivariant.

Assume the case
that $A$ is irreducible as a connection. Then the stabilizer
of $(A,0)$ is $\{\pm 1\}$ and this acts on the
$L^{2}_{2}(S\otimes_{\qu}E)$ factor only, by multiplication. 
Let $b\in X^{r}_{A}$ and set 
\begin{displaymath}
\pi_{A+a,\mPhi+\phi}' = \beta(a,\phi)(b,0) +\delta^{0}_{A,\mPhi}\xi(a,\phi; \beta(a,\phi)(b,0)).
\end{displaymath}
Then in the same manner as Lemma~\ref{pert1-lem} $\pi'$ is admissible provided the
support of $\beta$ is small, and by construction
$\pi_{A,0}' = (b,0)$. This defines perturbations in the connection
irreducible portion $\A^{*}$ of the reducible strata 
$\A\subset\CC$. 

Let us now consider the normal direction linearization $(L^{\nu}\pi)_{A}$ of
any perturbation $\pi$ at $A\in\A$. 
In preparation for  this we need a little
technical result:

\begin{lemma}\label{little-lem}
Let $V\to Y$ be a trivial 
real vector bundle of rank $\ge 2$ and let
$L\colon L^{2}_{2}(V)\to L^{2}_{2}(V)$ be a bounded linear operator. 
Regard $L^{2}_{2}(V)\subset C^{0}(V)$. Suppose that
$L(\sigma)(x)\in\langle \sigma(x)\rangle$ wherever $\sigma(x)\neq 0$. Then
there exists a real function $f\in L^{2}_{2}(Y)$ such that $L(\sigma)=f\sigma$
for all $\sigma$.
\end{lemma}

\proof
Let $\sigma$ be a nowhere zero section. Then $L(\sigma)=f\sigma$ for some $f\in C^{0}(Y)$.
Let $\sigma_{1}$ be a section which is pointwise linearly independent to $\sigma$ wherever
it is non-zero. Then $L(\sigma_{1})(x)=f_{1}(x)\sigma_{1}(x)$ for some $f_{1}$ at such
points. However it must also be the case that
\begin{displaymath}
L(\sigma+\sigma_{1})(x) = h(x) (\sigma+\sigma_{1})(x)
\end{displaymath}
for some $h$. If $\sigma_{1}(x)\neq 0$ this leads to the relation
\begin{displaymath}
(f(x) - h(x))\sigma(x) = (h(x) -f_{1}(x))\sigma_{1}(x)
\end{displaymath}
which forces $f(x) =h(x) = f_{1}(x)$. On the other hand if $\sigma_{1}(x) = 0$ then
we obtain
\begin{displaymath}
L(\sigma_{1})(x) = (h(x) -f(x))\sigma(x).
\end{displaymath}
Since we have the freedom to make other choices for $\sigma$ the only possibility
is that $h(x) =f(x)$ and so $L(\sigma_{1})(x) = 0$ wherever $\sigma_{1}(x)=0$.
Thus $L(\sigma_{1})=f\sigma_{1}$ i.e., $L(\sigma_{1})(x)=f(x)\sigma_{1}(x)$ 
\textit{for all} $x$.

Choose $\sigma_{1}$ to be nowhere vanishing and reverse the roles of $\sigma$
and $\sigma_{1}$ above. Then we obtain $L(k\sigma)=f(k\sigma)$ for any 
function $k\in L^{2}_{2}(Y)$. Finally, given \textit{any} section $\sigma'$
we may write this as a sum $\sigma'=k\sigma+\sigma_{1}$ where $\sigma$ and
$\sigma_{1}$ are as in the preceding paragraph. Then
\begin{displaymath}
L(\sigma') = L(k\sigma+\sigma_{1}) = f(k\sigma) +f\sigma_{1} =f\sigma'.
\end{displaymath}
If $L(\sigma)=f\sigma\in L^{2}_{2}$ for all $\sigma\in L^{2}_{2}$ then it must be the
case that $f\in L^{2}_{2}$ as well. \endproof

The next results limits the possibilities for the normal 
linearization of a perturbation which in turn forces it to be 
self-adjoint:

\begin{lemma}\label{linear-lem}
Given any admissible perturbation $\pi$ and $A\in \A$ there is a real function
$f_{A}\in L^{2}_{2}$ on {\it\,Y} such that
$(L^{\nu}\pi)_{A}(\delta\phi) =f_{A}\delta\phi$ for all 
$\delta\phi\in L^{2}_{2}(S\otimes_{\qu}E)$.
It follows that $(L^{\nu}\pi)_{A}$ is $L^{2}$--self-adjoint
on $L^{2}_{2}\subset L^{2}(S\otimes_{\qu}E)$.
\end{lemma}

\proof
Assume an admissible
perturbation $\pi$ is given. Then 
\begin{displaymath}
\langle (d_{A}\gamma,-\gamma(\phi)),\pi_{A,\phi}\rangle_{L^{2}} =0\quad
\mbox{\rm for all $\gamma\in L^{2}_{3}( \ad E)$}.
\end{displaymath}
Performing a variation $\phi\mapsto \phi+\delta\phi$ at $\phi=0$ gives
$\langle \gamma(\delta\phi),(L^{\nu}\pi)_{A}(\delta\phi)
\rangle_{L^{2}}=0$ for all $\gamma$. Write $(L^{\nu}\pi)_{A}(\delta\phi)=\delta\psi$.
By assumption $\delta\phi,\delta\psi\in L^{2}_{2}\subset C^{0}$ so we can consider them as
continuous sections.
Then pointwise we have $\langle \gamma(\delta\phi),\delta\psi\rangle_{x}=0$.
A local model for the fibre of $S\otimes_{\qu}E$ is just $\qu$ and with
the action of $\gamma$ as multiplication by ${\rm Im\,}\qu$. 
Thus we see that $\delta\psi(x)=(L^{\nu}\pi)_{A}(\delta\phi)(x)\in\langle\delta\phi(x)\rangle$ 
at all points $x$ where $\delta\phi(x)\neq 0$. The proof is completed by
Lemma~\ref{little-lem}. \endproof

Let us now construct perturbations normal to $\A\subset\CC$. Assume the cutoff
$\beta$ on $X_{A,0}$ is invariant under the stabilizer action. 
Let $f_{A}$ be a real $L^{2}_{2}$ function on $Y$. Set
\begin{eqnarray*}
\pi_{A+a,\mPhi+\phi}'' &=& \beta(a,\phi)(0,f_{A}\phi) 
+\delta^{0}_{A,\mPhi}\xi(a,\phi;\beta(a,\phi)(0,f_{A}\phi))\\
&=&\beta(a,\phi)f_{A}\phi\in X_{A+a,\mPhi+\phi}.
\end{eqnarray*}
This is
again an admissible perturbation for  $\supp\beta$ small and 
the linearization of $\pi''$ in a normal
direction $\delta\phi$ at $(A,0)$ is
$(L^{\nu}\pi'')_{A}(\delta\phi) = f_{A}\delta\phi$.

Thus we have:

\begin{prop}\label{pert-2}
If $A$ is irreducible then 
given any $b\in X^{r}_{A}$ there is an
admissible perturbation $\pi'$ such that 
$\pi_{A,0}'=(b,0)\in X_{A,0}=X^{r}_{A}\times L^{2}_{2}(S\otimes_{\qu}E)$.
On the other hand for any $A$ there exists
an admissible perturbation $\pi''$ such that 
$\overline{\pi}''=0$ and
$(L^{\nu}\pi'')_{A}(\delta\phi)=f_{A}\delta\phi$
given any real function $f_{A}\in L^{2}_{2}(Y)$.
Furthermore the support of $\pi'$ and $\pi''$ may be chosen
to be contained in an arbitarily small $\G$--invariant 
neighbourhood of the orbit  $\G\cdot(A,0)$ in $\CC$.
\end{prop}

\proofof{Proposition~\ref{adm2-lem}}$\phantom{99}$\nl
Let 
$\mathcal{X}^{r}_{\overline{\pi}}(A)=*F_{A}+*\overline{\pi}_{A}\in X^{r}_{A}$. Then 
$(\mathcal{X}^{r}_{\overline{\pi}})^{-1}(0)/\G_{E}=\M^{r}_{\overline{\pi}}$. 
Let $H^{\tau}_{A}$ denote the cokernel of 
$(L\mathcal{X}^{r}_{\pi})_{A}\colon L^{2}_{2}(\mLambda^{1}\otimes\ad E)\to X^{r}_{A}\cap L^{2}_{1}$. Then
at the reducible $A=(A,0)$  the cokernel $H^{2}_{A}$  of 
$(L\mathcal{X}_{\pi})_{A,0}$
splits as $H^{\tau}_{A}\oplus H^{\nu}_{A}$. 
$H^{\nu}_{A}$ is the cokernel of the normal operator $D^{\pi}_{A}$.
(Recall the map $\mathcal{X}_{\pi}=\mathcal{X} +\pi$
of (\ref{chi-map}) and its linearization (\ref{chi-lin}).)

\textbf{Step 1}\qua 
For a ZHS the orbit of the trivial connection $[\mTheta]\in\M^{r}$ is already
isolated in $\B_{\A}$ since $H^{\tau}_{\mTheta}\cong H^{1}(Y)=\{0\}$.
By Proposition~\ref{pert-2} and the compactness of $\M^{r*}$, we can find a finite set of perturbations
$\{ \pi^{(i)}\}$ with support away from $\{[\mTheta]\}$ such that if $v\in H^{\tau}_{A}$, 
$A\in (\mathcal{X}^{r})^{-1}(0)\cap\A^{*}$,  is $L^{2}$--orthogonal
to each $\pi^{(i)}_{A}$ then $v=0$. Thus by Sard--Smale there is
a perturbation, call it $\pi_{1}$ so that $(\mathcal{X}^{r}_{\overline{\pi}_{1}})^{-1}(0)$ is
cut out equivariantly transversely over $\A^{*}$, i.e.\ $H^{\tau}_{A}=\{0\}$ for every
$A\in (\mathcal{X}^{r}_{\overline{\pi}_{1}})^{-1}(0)\cap \A^{*}$. Hence $\M^{r}_{\overline{\pi}_{1}}$
is, by the local Kuranishi model, a finite set of points which are 
\textit{non-degenerate within $\B_{\A}$}. 

\textbf{Step 2}\qua
Let $[A]\in\M^{r}_{\overline{\pi}_{1}}$. The normal operator $D^{\pi_{1}}_{A}$ at $A$
is of the form $D_{A}+f_{A}$, by Lemma~\ref{linear-lem}. This operator is self-adjoint
Fredholm and therefore has discrete spectrum. Let $\pi_{2}$ be a perturbation
with the property that $\overline{\pi}_{2}=0$ and 
$(L^{\nu}\pi_{2})_{A}\delta\phi=\mu_{A}\delta\phi$,
$\mu_{A}\in\R$ where $|\mu_{A}|$ is less than the distance of the closest non-zero
eigenvalue of $D^{\pi_{1}}_{A}$ from zero. Then $D^{\pi_{1}+\pi_{2}}_{A}$ has trivial
kernel and $[A]$ is a non-degenerate point in $\M_{\pi_{1}+\pi_{2}}$. $\pi_{2}$
can be chosen to have support in an arbitarily small $\G$--invariant neighbourhood
of the orbit of $A$. Repeating this procedure for every $[A]$ we 
can find a perturbation $\pi'$ such that
$\M^{r}_{\overline{\pi}'}$ consists entirely of non-degenerate
points within $\B$.

\textbf{Step 3}\qua
After the preceding steps, $\M^{r}_{\overline{\pi}'}$ is
isolated in $\M_{\pi'}$. By Proposition~\ref{pert-1} 
and the compactness of $\M^{*}_{\pi'}$ we can find
a finite set of perturbations $\{\pi^{(j)}\}$ supported away from
$\M^{r}_{\overline{\pi}'}$ such that if $v\in H^{2}_{x}$,
$x\in\mathcal{X}_{\pi'}^{-1}(0)\cap\CC^{*}$, is
$L^{2}$--orthogonal to every $\pi^{(j)}$, then $v=0$. Thus by Sard--Smale
there exists a $\pi_{1}'$  which is an arbitarily small linear combination of these
$\pi^{(j)}$'s, such that $\mathcal{X}_{\pi_{1}+\pi_{1}'}^{-1}(0)$ is cutout equivariantly
transversely over $\CC^{*}$, i.e.\ $H^{2}_{x}=\{0\}$ for every 
$x\in\mathcal{X}_{\pi_{1}+\pi_{1}'}^{-1}(0)\cap\CC^{*}$. Note that
$\pi_{1}'$ is supported away from $\M^{r}_{\overline{\pi}'}$. 
Choosing our final perturbation $\pi$ to be $\pi'+\pi_{1}'$ we get
$\M_{\pi}$ non-degenerate. 

At every stage in Steps 1, 2 and 3 we can make the chosen perturbation
as small as we like in the uniform norm 
\begin{displaymath}
\normi{\pi'}{\B}
= \sup_{A,\mPhi}\Bigl\{ \normi{\pi_{A,\mPhi}'}{L^{2}_{2,A}}\Bigr\}.
\end{displaymath}
This completes the proof
of Proposition~\ref{adm2-lem}.

\section{Proof of Theorem~\protect\ref{thrm1}}\label{proof1-sect}

Let $(g_{0},\pi_{0})$ and $(g_{1},\pi_{1})$ be given. Assume that
$\pi_{i}$ is non-degenerate with respect to $g_{i}$. In order
to compare the moduli spaces for different metrics we may assume,
as in \S\ref{comp-sect} a fixed model for the spinor bundle with respect to
$g_{0}$. Then we have a SW--equation depending smoothly on the parameter $t$
corresponding to the metric $g_{t}=(1-t)g_{0}+t g_{1}$ and with
perturbation $\pi_{t}=(1-t)\pi_{0}+t\pi_{1}$. In this section we
shall assume objects sub- or superscripted with \lq$t$\rq\ are
with respect to $g_{t}$. 

To the family $\{(g_{t},\pi_{t})\}$ we have
a {\it parameterized moduli space}
\begin{displaymath}
Z = 
\bigcup_{t}\M_{g_{t},\pi_{t}}\times\{t\}\subset\CC\times[0,1].
\end{displaymath}
As in \cite{BH} and \cite{lim1} to prove invariance of $\tau(Y)$ we need to show that $Z$
is, after suitable perturbation, a compact 1--dimensional cobordism with the appropriate
singularities. The counter-terms in the definition of $\tau(Y)$ are due
to these singularities. 

In our analysis of $Z$ we work first
with the reducible strata $Z^{r}$. In the following the notation
${Z}^{r*}$ denotes the connection irreducible portion of ${Z}^{r}$.
In the parameterized context
an {\it admissible time-dependent perturbation} $\sigma$ is one which is a finite
sum $\sum_{i}\varrho_{i}(t)\pi^{(i)}$ where $\pi^{(i)}$ is admissible
and $\varrho_{i}$ has support in $[0,1]$.  
${Z}_{\sigma}$, ${Z}^{r}_{\sigma}$ etc. shall denote
perturbed parameterized moduli spaces.
Recall the uniform norm, for non-time-dependent perturbations,
$$
\normi{\pi}{\B}
= \sup_{A,\mPhi}\Bigl\{ \normi{\pi_{A,\mPhi}}{L^{2}_{2,A}}\Bigr\}.
$$

\begin{lemma}\label{inv-lem1}
There exists  an admissible time-dependent perturbation $\sigma$
such that the perturbed parameterized reducible moduli space 
$Z^{r}_{\sigma}$ is non-degenerate as a subspace of $\B_{\A}\times [0,1]$. 
Furthermore if $\normi{\pi_{0,1}}{\B}<\delta_{0}$ 
then we can assume 
$\normi{\pi_{t}+\sigma(t)}{\B}<2\delta_{0}$.
\end{lemma}

$Z^{r}_{\sigma}$ can be regarded as the $\G$--quotient of the
zeros of the map $\mathcal{X}^{r}_{\sigma}(A,t)=*_{t}F_{A}+\overline{\sigma}(t)_{A}\in X^{r}_{A}$.
The proof of the Lemma 
follows easily from constructing and applying time-dependent perturbations
to $\mathcal{X}^{r}_{\sigma}$ supported away from $\{[\mTheta]\}$
in the manner of \S\ref{pert-sect}.
In this way the strata corresponding to the trivial
connection is isolated in $Z^{r}_{\sigma}$ and
 is the product $\{[\mTheta]\}\times [0,1]$. The
irreducible portion of $Z^{r}_{\sigma}$ with the choice of
$\sigma$ in Lemma~\ref{inv-lem1} is a compact corbordism
between $\M^{r*}_{g_{0},\pi_{0}}$ and $\M^{r*}_{g_{1},\pi_{1}}$.

Assume now ${Z}^{r}_{\sigma}$ as in Lemma~\ref{inv-lem1}.
The normal operator at $A$ with respect to $(g_{t},\pi_{t}+\sigma(t))$ will be
denoted by $D_{A}^{t,\sigma}$ (Eq.~(\ref{split-eq})). The kernel of this operator ($=$ cokernel
by Lemma~\ref{linear-lem}) is the \textit{normal cohomology} $H^{\nu}_{A,t}$.

Let $u\mapsto ([A(u)],t(u))$, $|u|<\epsilon$ be a 1--1 parameterization
of an open subset ${\cal J}$ of ${Z}^{r}_{\sigma}$. 
Let ${\bf K}=\R$ if ${\cal J}$ is in the connection-irreducible strata and
${\bf K}=\qu$ if ${\cal J}$ is in the connection-trivial strata.

\begin{defn}\rm
Call ${Z}^{r}_{\sigma}$ {\it normally transverse}
along ${\cal J}$ if the family $\{ D^{t(u),\sigma}_{A(u)}\}$
has transverse spectral flow as ${\bf K}$--linear operators. 
(Recall that transverse spectral flow is the situation of simple
eigenvalues, with respect to $\bf K$, crossing
zero transversely.)
Call ${Z}^{r}_{\sigma}$ {\it normally tranverse}
if it is normally transverse in a neighbourhood of every point.
\end{defn}

In terms of local models, let
$A=A(0)$, $t_{0}=t(0)$ and  $\mathcal{U}$ be a 
sufficiently small $\stab(A_{0})$--invariant neighbourhood of $(A_{0},t_{0})$
in the slice  $(A_{0},t_{0})+X^{r}_{A}\times[0,1]$. $\mathcal{U}$ can be identified with
a neighbourhood of
$([A_{0}],t_{0})$ in $\B_{A}\times[0,1]$. 
To simplify notation henceforth denote $H^{\nu}_{A_{0},t_{0}}$
by $H^{\nu}_{0}$.
Assume
$H^{\nu}_{0}$ is non-trivial, otherwise $\mathcal{U}$
can be chosen such that $D^{t,\sigma}_{A}$ is invertible in $\mathcal{U}$. 
Consider the restriction of $D^{t,\sigma}_{A}$
to the normal cohomology $H^{\nu}_{0}$ followed by $L^{2}$--projection
$\prod$ back onto $H^{\nu}_{0}$. This determines, for each $(A,t)\in\mathcal{U}$
a symmetric operator $T(A,t)$  acting on 
$H^{\nu}_{0}$. The latter space is endowed with the natural
$L^{2}$--inner product. Then the kernel (cokernel) of $D^{t,\sigma}_{A}$
is exactly modelled by the kernel (cokernel)  of $T(A,t)$. 
Denote the symmetric operators which commute with $\mathbf{K}$ by
$\mathrm{Sym}_{{\bf K}}(H^{\nu}_{0})$. 
Let $u\mapsto (A(u),t(u))\in X^{r}_{A}\times(0,1)$, $|u|<\epsilon$ be a 1--1 parameterization
of an open subset $\mathcal{J}=\mathcal{U}\cap (\mathcal{X}^{r}_{\sigma})^{-1}(0)$
of ${Z}^{r}_{\sigma}$.
Then the
condition of being normally transverse along 
$\mathcal{J}$ translates as (i) $H^{\nu}_{0}\cong{\bf K}$ and (ii) the path 
$u\mapsto T(A(u),t(u))\in{\rm Sym}_{{\bf K}}(H^{\nu}_{0})\cong\R$
is transverse to $\{0\}$.

\begin{lemma}\label{inv-lem2}
Assume $\sigma$ as in Lemma~\ref{inv-lem1}.
There exist an admissible time-dependent perturbation $\sigma'$ such that 
${Z}^{r}_{\sigma+\sigma'}\simeq{Z}^{r}_{\sigma}$ and
${Z}^{r}_{\sigma+\sigma'}$ is normally transverse.
Furthermore if $\normi{\pi_{0,1}}{\B}<\delta_{0}$ we can assume 
$\normi{\pi_{t}+\sigma(t)+\sigma'(t)}{\B}<3\delta_{0}$.
\end{lemma}

\proof
We divide the argument into the separate cases of the irreducible and trivial
strata of ${Z}^{r}_{\sigma}$. No matter what perturbation $\sigma$ is
chosen the trivial strata is always $\{[\mTheta]\}\times[0,1]$. However changing
$\sigma$ can change ${Z}^{r*}_{\sigma}$.
The space $\mathcal{S}$  of admissible 
perturbations is a normed linear space, with the norm
$\normi{\cdot}{\B}$.
Since ${Z}^{r*}_{\sigma}$ is already non-degenerate as a subspace
of $\B^{*}_{\A}\times [0,1]$
i.e.\ $(\mathcal{X}^{r}_{\sigma})^{-1}(0)\cap\A^{*}\times[0,1]$ 
is cutout equivariantly transversely, it follows that for any
sufficiently uniformly small $\sigma'\in\mathcal{S}$, ${Z}^{r*}_{\sigma}$ and 
${Z}^{r*}_{\sigma+\sigma'}$ are related by a cobordism which is
a product and thus  are diffeomorphic spaces. In fact the transverse
condition means that the normal bundle to ${Z}^{r*}_{\sigma}$ in $\B_{\A}\times [0,1]$
at any point is isomorphic to $\mathcal{S}/\mathcal{S}_{0}$ where $\mathcal{S}_{0}$
is the subspace of those $\pi$ such that $\overline{\pi}=0$.

\textbf{Case 1: Irreducible strata\ }
Fix $([A_{0}],t_{0})\in{Z}^{r*}_{\sigma}$, $t\neq 0,1$ and 
let $\mathcal{U}$ be a sufficiently small $\G$--invariant neighbourhood of
$(A_{0},t_{0})$ in $(A_{0},t_{0})+X^{r}_{A}\times[0,1]$ 
such that a local model for
$D^{t,\sigma}_{A}$ as above exists in $\mathcal{U}$.
Assume $H^{\nu}_{0}$ is non-trivial. We examine the effect of a perturbation
on the family $D^{t,\sigma}_{A}$ along ${Z}^{r*}_{\sigma}$.  

Consider the parameterized local model map $P\colon \mathcal{U}\times \mathcal{S}\to
{\rm Sym}_{\R}(H^{\nu}_{0})$  based at $(A_{0},t_{0},0)$ given by
$$
P(A,t,\pi)=
\textstyle{\prod}\circ D^{t,\sigma+\beta(t)\pi}_{A}
$$
where $\beta=\beta(t)$ is a cutoff function on $\R$ with support close to $t_{0}$.
Note that $D^{t,\sigma+\pi}_{A}=D^{t}_{A}+(L^{\nu}\sigma(t)+L^{\nu}\pi)_{A}$.
By Lemma~\ref{linear-lem}, $(L^{\nu}\sigma(t))_{A}\phi=f_{A,t}\phi$ and
$(L^{\nu}\pi)_{A}\phi=h_{A}\phi$ for some functions $f_{A,t}$ and 
$h_{A}$ on $Y$.
Since $\mathcal{S}$ is a linear space, we may identify tangent vectors 
$\delta\pi$ with elements $\pi$ in $\mathcal{S}$.
We then have, for the derivative of $P$ at $(A_{0},t_{0},0)$,
$$
dP(\delta\pi)\phi = \textstyle{\prod}\Bigl( h_{A_{0}} \phi\Bigr),\quad
dP(\delta a)\phi= \textstyle{\prod}\Bigl( \delta a\cdot\phi+ (Lf)_{A_{0},t_{0}}(\delta a) \phi\Bigr).
$$
By choosing $h=h_{A_{0}}=1$ we see that the image of $dP$ includes
at least the span of the identity operator in ${\rm Sym}_{\R}(H^{\nu}_{0})$; thus
$\mbox{\rm rank}(dP)\ge 1$.

\textbf{Claim}
\textit{If $\dimn_{\R}(H^{\nu}_{0})>1$ then $\mbox{\rm rank}(dP)\ge 2$}.

In order to establish the claim we invoke the unique 
continuation principle for $H^{\nu}_{0}$ i.e.\ if $\phi\in H^{\nu}_{0}$
then $\phi$ cannot vanish on an open set unless $\phi=0$. Writing 
$A_{0}=\mTheta +a$ where $\mTheta$ is smooth, then $\phi\in 
H^{\nu}_{0}$ is a solution of the perturbed smooth Dirac operator
$D_{\mTheta}\phi+a\cdot\phi +f\phi=0$ where $a$ and $f$ are continuous. 
Unique continuation holds for such solutions.

Let $\{\phi_{1},\dots,\phi_{n}\}$,
$n>1$
be a $\R$--orthonormal basis for $H^{\nu}_{0}$. The matrix
of $dP(\delta\pi)$ with respect to this basis is 
$( \langle h\phi_{i},\phi_{j}\rangle_{L^{2}})$. Assume the
rank of $dP$ is unity. This implies that
$\langle h\phi_{i},\phi_{j}\rangle_{L^{2}}=0$ for all $h$
and $i\neq j$. 
This in turn implies the pointwise orthogonal condition
$\langle{\phi_{i},\phi_{j}}\rangle_{y}=0$, $i\neq j$ for all $y\in Y$.
It then follows that $\langle dP(\delta a)\phi_{i},\phi_{j}\rangle_{L^{2}} = 
\langle \delta a\cdot\phi_{i},\phi_{j}\rangle_{L^{2}}$, $i\neq j$.
However the Clifford
action of $\mLambda^{1}\otimes  \ad E$ on $S\otimes_{\qu} E$ 
is fibrewise  transitive. Thus we can find a $\delta a$ such that
$\langle \delta a\cdot\phi_{i},\phi_{j}\rangle_{L^{2}}\neq 0$, $i\neq j$. 
This proves that the image of $dP$ is not contained in the span
of the identity in ${\rm Sym}_{\R}(H^{\nu}_{0})$. Therefore $dP$ is at least
rank two and the claim is proven.

Let $u\mapsto (A(u),t(u))\in X^{r}_{A}\times (0,1)$, $|u|<\epsilon$ be  a 1--1 parameterization of 
$\mathcal{U}\cap (\mathcal{X}^{r}_{\sigma})^{-1}(0)$ with $A(0)=A_{0}$ and $t(0)=t_{0}$.
Let $T_{\sigma}(u)=P(A(u),t(u),0)$ be the local model for $D^{t(u),\sigma}_{A(u)}$ on
${\rm Sym}_{\R}(H^{\nu}_{0})$. 
By construction, $P(A_{0},t_{0},0)$ is the zero
operator on ${\rm Sym}_{\R}(H^{\nu}_{0})$. In ${\rm Sym}_{\R}(H^{\nu}_{0})$
the space of invertible operators is a codimension one real variety ${\cal V}$.
Any point which is not the zero operator in this variety represents
an operator of non-trivial rank. 

Let $X^{r}$ be the vector bundle over $\A^{*}$ whose fiber at $A$ is the slice
space $X^{r}_{A}$ and let $Q\colon \mathcal{U}\times \mathcal{S}\to
X^{r}\times{\rm Sym}_{\R}(H^{\nu}_{0})$ be given by
$$
Q(A,t,\pi) = (\mathcal{X}^{r}_{\sigma+\beta(t)\pi}(A), P(A,t,\pi) ).
$$
Since this is a submersion onto the first factor along
$(\mathcal{X}^{r}_{\sigma})^{-1}(0)\cap\A^{*}\times[0,1]$ (the transversality
condition) and $\mbox{\rm rank}(dP)\ge 2$ if $\dimn_{\R}(H^{\nu}_{0})>1$ and is onto if
$\dimn_{\R}(H^{\nu}_{0})=1$, then there is a time-dependent 
perturbation $\sigma_{0}(t)\colon=\beta(t) \pi$
such that the deformation of the family 
$\{T_{\sigma+s\sigma_{0}}\}$ at $s=0$
is normal to the path $T_{\sigma}=T_{\sigma}(u)$. 
Therefore we can choose an arbitarily small $\sigma'$ so that the
operators $T_{\sigma+\sigma'}(u)$ have non-trivial rank for
all $u$. (Note: at this stage we do not have sufficently many perturbations in hand to
make $T_{\sigma+\sigma'}$ transverse to ${\cal V}$.)  Thus if we work with $\sigma+\sigma'$
we find that the rank (over $\R$) of $H^{\nu}_{A(u),u}$ near $u=0$ drops
by one if $\dimn_{\R}(H^{\nu}_{0})>1$ and becomes transverse to ${\cal V}=\{0\}$
if $\dimn_{\R}(H^{\nu}_{0})=1$. 
To complete the argument to obtain normal transversality
globally over the connection-irreducible strata, proceed by an induction
argument with the overall rank of $H^{\nu}_{A,t}$ over
${Z}^{r*}_{\sigma+\sigma'}$ decreasing by one in each step.
Letting $\sigma'$ denote the final perturbation we see that over
${Z}^{r*}_{\sigma+\sigma'}$ there exists a finite number of
points where $H^{\nu}_{A,t}$ is non-trivial and these points $H^{\nu}_{A,t}\cong\R$
and with $T_{\sigma+\sigma'}$ transverse to ${\cal V}=\{0\}$. This
is equivalent to transverse spectral flow.
The last assertion of the lemma in this case is a consequence of the observation
that the induction is completed in a finite number of steps
and in each step we may take the perturbation to be as small
as we like. 

\textbf{Case 2: Trivial strata\ } Let $([A_{0}],t_{0})\in \{[\mTheta]\}\times[0,1]$.
Here the relevant parameterized local model map $P$ is the same as the map $P$
as above but with $A=\mTheta$ fixed, i.e.\ $P\colon\mathcal{S}\to 
{\rm Sym}_{\qu}(H^{\nu}_{0})$. The argument proceeds just as before 
(but without the complication of the deformation in the moduli space)
provided
we can again establish that if  $\dimn_{\qu}(H^{\nu}_{0})>1$ then
$\mbox{\rm rank}(dP)\ge 2$. This time let
$\{\phi_{1},\dots,\phi_{n}\}$, $n>1$
be a $\qu$--orthonormal basis for $H^{\nu}_{0}$. Again if we assume
the rank of $dP$ is unity we get the pointwise orthogonal condition
$\langle{\phi_{i},\phi_{j}}\rangle_{y}=0$, $i\neq j$ for all $y\in Y$.
However this would mean that $S\otimes_{\qu}E$ has at least 8
pointwise orthogonal non-zero sections. This is impossible since
$S\otimes_{\qu}E$ is rank 4. 

This completes the proof of the lemma. \endproof

\begin{remark}\rm
A more satisfactory result would be that $P$ is a submersion onto
${\rm Sym}_{{\bf K}}(H^{\nu}_{0})$ which is the situation in \cite{BH};
then transverse spectral flow follows easily by Sard--Smale.
A submersion does not seem to be generally true in our and the original SW context.
The same problem is encountered in \cite{lim2} and \cite{lim3}.
\end{remark}

\begin{defn}\rm
Suppose ${Z}^{r}_{\sigma}$ is normally tranverse and let
$u\mapsto ([A(u)],t(u))$, $|u|<\epsilon$ be a 1--1--parameterization
of an open neighbourhood in ${Z}^{r}_{\sigma}$.
A point in $Z^{r}_{\sigma}$ which is contained in such a 
parameterization and 
where there is spectral flow for $D^{t(u),\sigma}_{A(u)}$
is called a {\it singular} or {\it bifurcation } point. 
\end{defn}

At a singular 
point $([A_{0}],t_{0})$, the local model for 
${Z}_{\sigma}$ is the quotient by $\stab (A_{0})$ of 
the zeros of a $\stab (A_{0})$--equivariant obstruction map
$\Xi\colon H^{\nu}_{0}\times\R
\to H^{\nu}_{0}$ of the
form
$$
\Xi(q,t) = qt.
$$
(See \cite{lim2} and \cite{BH}.)
This in turn implies that the a neighbourhood of $([A_{0}], t_{0})$
is the zeros of the map $[0,\infty)\times\R\to\R$, $(r,t)\mapsto rt$
with $\{ 0\}\times \R$ corresponding to the reducible portion and
$(0,\infty)\times\{0\}$ the irreducible. One other consequence of the
local model in this normal transverse situation is that the points
corresponding to the irreducibles sufficiently near $([A_{0}],t_{0})$
are non-degenerate.

On the other hand, at a non-singular point $([A_{0}],t_{0})$
of a normally transverse ${Z}^{r}_{\sigma}$
the Kuranishi local model gives a neighbourhood of $([A_{0}],t_{0})$
in $\B\times [0,1]$ an isolated open interval.

\begin{cor}\label{inv-cor1}
Assume $\sigma$ as in Lemma~\ref{inv-lem1}.
There exists a time-dependent admissible perturbation $\sigma'$ such that 
(i) ${Z}^{r}_{\sigma+\sigma'}\simeq{Z}^{r}_{\sigma}$
(ii) ${Z}^{r}_{\sigma+\sigma'}$ is normally transverse and
(iii) ${Z}^{*}_{\sigma+\sigma'}$ is non-degenerate.
Furthermore if $\normi{\pi_{0,1}}{\B}<\delta_{0}$ we can assume 
$\normi{\pi_{t}+\sigma(t)+\sigma'(t)}{\B}<4\delta_{0}$. 
\end{cor}

\proof
Run through the proof of Lemma~\ref{inv-lem2}. The comments above
tell us that ${Z}^{*}_{\sigma+\sigma'}$ is non-degenerate in
a neighbourhood of ${Z}^{r}_{\sigma+\sigma'}$. Now
construct and apply admissible time-dependent perturbations 
$\sigma''$ in the manner of \S\ref{pert-sect}, which can be chosen to
have support away from ${Z}^{r}_{\sigma+\sigma'}$, making
all of ${Z}^{*}_{\sigma+\sigma'+\sigma''}$ non-degenerate.
The perturbation $\sigma''$ can be chosen arbitarily small. \endproof

\textbf{Completion of proof of Theorem~\ref{thrm1}} 
As above we have two non-degenerate metrics and perturbations
$(g_{0},\pi_{0})$ and $(g_{1},\pi_{1})$ where $\pi_{i}$ is small
with respect to $g_{i}$.

Assume first the case that the metric $g=g_{0}=g_{1}$ is unchanging. The condition
$\pi_{0}$, $\pi_{1}$ are \textit{small} (Definition~\ref{small-defn}) implies
$\M^{r*}_{\pi_{i}}\subset\cup_{j}\N^{j}$ where the $\N^{j}$ are
as in the definition of the proposed invariant. By 
Corollary~\ref{inv-cor1} we can find a parameterized
moduli space ${Z}_{\sigma}$ such that 
\begin{enumerate}
\item ${Z}^{r*}_{\sigma}$
is a smooth compact 1--dimensional corbodism between $\M^{r*}_{\pi_{0}}$ and
$\M^{r*}_{\pi_{1}}$. Additionally we know from \cite{taubes} that this is
an oriented cobordism so that it's boundary is $\M^{r*}_{\pi_{1}}-\M^{r*}_{\pi_{0}}$
where $\M^{r*}_{\pi_{0,1}}$ are given Taubes' orientation
\item ${Z}^{r*}_{\sigma}\subset \dis \cup_{j}\N^{j}$
\item $\overline{{Z}^{*}_{\sigma}}$ is a smooth compact 1--manifold
with boundary $$\M^{*}_{\pi_{0}}\cup \M^{*}_{\pi_{1}}\cup\{
\textrm{singular points in $Z^{r}_{\sigma}$}\}.$$
\end{enumerate} 

Just as in \cite{BH} it is seen that
\begin{eqnarray}
\lefteqn{\sum_{[A]\in\M^{r*}_{\pi_{1}}}\epsilon[A]\kappa[A]
-\sum_{[A]\in\M^{r*}_{\pi_{0}}}\epsilon[A]\kappa[A]} \label{inv-equ1}\\
&=&\#\{\textrm{singular points on ${Z}^{r*}_{\sigma}$}\}\bmod 2.\nonumber
\end{eqnarray}
For completeness we give an argument.
Fix a component $\N^{j}$ and consider ${Z}^{r*}_{\sigma}\cap \N^{j}$.
In the definition of $\kappa[A]$ for $[A]\in \M^{r*}_{\pi_{0,1}}\cap\N^{j}$
choose all the paths $[\gamma]$ to be in the same homotopy class rel 
$\{[\mTheta]\}\cup \N^{j}$.
Then for these $[A]$'s the term 
$\cs(\overline{\gamma}(0))-\cs(\overline{\gamma}(1))$ is the same constant.
Make this choice.
Then $\kappa[A]$ is 
the normal spectral flow $\SF^{\nu}(\gamma)$ from
$[\mTheta]$ to $[A]$ in the given fixed homotopy class of $[\gamma]$ plus
a fixed additive constant. Notice then that $\kappa[A]$ changes exactly by
the normal spectral flow as we vary $[A]$ within $\N^{j}$.
Let $\Gamma$ be a connected
component of ${Z}^{r*}_{\sigma}\cap \N^{j}$
with non-empty boundary 
$\{[A],[A']\}\subset \M^{r*}_{\pi_{0,1}}\cap\N^{j} $.
After some consideration it is seen that the three following 
sums compute the $\bmod 2$ normal spectral flow along $\Gamma$
and thus the $\bmod 2$ cardinality of the singular points on $\Gamma$:
\begin{enumerate}
\item $\epsilon[A]\kappa[A]+\epsilon[A']\kappa[A']=\pm(\kappa[A]-\kappa[A'])$ 
when $[A], [A']\in \M^{r*}_{\pi_{1}}$
\item $-\epsilon[A]\kappa[A]-\epsilon[A']\kappa[A']=\pm(\kappa[A]-\kappa[A'])$ 
when $[A], [A']\in \M^{r*}_{\pi_{0}}$
\item $\epsilon[A]\kappa[A]-\epsilon[A']\kappa[A']=\pm(\kappa[A]-\kappa[A'])$
when $[A]\in \M^{r*}_{\pi_{1}}$, $[A']\in \M^{r*}_{\pi_{0}}$.
\end{enumerate}
On the other hand, if $\Gamma$ has empty boundary then the number of singular points on $\Gamma$ equals
the normal spectral flow around $\Gamma$ and this is zero, since it is contained
within $\mathcal{N}^{j}$.
From this it is straightforward to deduce (\ref{inv-equ1}) by rearranging the sum.

Next we compute that the difference
\begin{eqnarray}
\lefteqn{ \frac{1}{2}c(g,\pi_{1})-\frac{1}{2}c(g,\pi_{0})}\label{inv-equ2}\\
&=&\textrm{$\qu$--spectral flow of 
$\{D^{t,\sigma}_{\mTheta} \}_{t=0}^{1}$}
\nonumber\\
&=&\pm\#\Bigl\{\textrm{singular points on trivial strata $\{[\mTheta]\}\times[0,1]$}\Bigr\}
\nonumber.
\end{eqnarray}
Finally we have equality of the sums
\begin{equation}\label{inv-equ3}
\sum_{\M^{r*}_{\pi_{1}}}\frac{1}{4}\epsilon[A] = 
\sum_{\M^{r*}_{\pi_{0}}}\frac{1}{4}\epsilon[A]
\end{equation}
both being $1/2$ of the algebraic sum  which is
Casson's invariant \cite{taubes}. Thus from (\ref{inv-equ1}), (\ref{inv-equ2}),
(\ref{inv-equ3}) we find that
\begin{eqnarray*}
&&\frac{1}{2}c(g,\pi_{1}) +
\sum_{\M^{r*}_{\overline{\pi}_{1}}} \epsilon[A]\Bigl(\kappa[A]+\frac{1}{4}\Bigr)\\
&&\mbox{}-\frac{1}{2}c(g,\pi_{0})
-\sum_{\M^{r*}_{\overline{\pi}_{0}}} \epsilon[A]\Bigl(\kappa[A]+\frac{1}{4}\Bigr)\\
&&\equiv  \#\Bigl\{\textrm{singular points on ${Z}^{r}_{\sigma}$}\Bigr\}\bmod 2\\
&&\equiv  \sum_{\M^{*}_{\pi_{0}}} 1 - \sum_{\M^{*}_{\pi_{1}}}1\quad \bmod 2 .
\end{eqnarray*}
The last line follows from $\overline{{Z}^{*}_{\sigma}}$ being a smooth compact 1--manifold
with boundary $\M^{*}_{\pi_{0}}\cup \M^{*}_{\pi_{1}}\cup\{
\mbox{\rm singular points on $Z^{r}_{\sigma}$}\}$.
Thus the independence of $\tau(Y)$ on choice of small, 
non-degenerate perturbation $\pi$
is established. 

The general case  $g_{0}\neq g_{1}$ follows
an identical argument except for the following details. 
When varying the metric spectral flow can occur at the trivial connection $\mTheta$ 
in $\SF^{\nu}(\gamma)$, which is the initial point of $\gamma$. 
However the operator $D_{\mTheta}$ at this point is quaternionic and 
thus there is no change $\bmod 2$. 
Secondly the neighbourhoods
$\N^{j}$ are defined with reference to the background metric, thus
we get for the different metrics $g_{0}$, $g_{1}$ two sets of neighbourhoods
$\N^{j}_{0}$, $\N^{j}_{1}$. However by what we have established we can
make any choice of (non-degenerate) $\pi_{0,1}$ we like.
Choose $\pi_{0,1}$ sufficiently small in norm so that 
$\M^{r*}_{\pi_{0,1}}\subset \cup_{j}(\N^{j}_{0}\cap\N^{j}_{1})$.
Then we may proceed with the rest of the argument as before. This
proves that $\tau(Y)$ is an invariant.

Finally, let us show that $\tau(-Y)=\tau(Y)$, $-Y$ denoting $Y$ with the
reversed orientation.
Reversing orientation but
keeping the metric, spin structure $P\to Y$ and and spinor bundle $S$
fixed simply changes the action of Clifford mutiplication by $-1$.
The SW--equation of the orientation reversed structure is the same as
the orginal except that the Dirac operator $D_{A}$ switches to $-D_{A}$.
If $\pi=(*k,l)$ is the non-degenerate and small perturbation used
to compute $\tau(Y)$ then choose $\pi'=(*k,-l)$ for the reversed structure.
Thus if $(A,\mPhi)$ is a SW--solution with respect to $\pi$ then
$(A,-\mPhi)$ is a solution of the orientation reversed
situation for $\pi'$.
In the following $\M^{-}$, $\epsilon^{-}$ etc. will refer to the reversed
orientation structure.
Thus $\M^{-}_{\pi'}=\M_{\pi}$ and $\pi'$ is a non-degenerate small
perturbation for $\M^{-}$.

The normal and tangential deformation operators $D^{\pi}_{A}$ and 
$K^{\overline{\pi}}_{A}$ in the reversed situation are
the negatives of those in the original. Then 
$\SF^{\nu,-}(\gamma)=-\SF^{\nu}(\gamma)-\dimn\kernel 
K^{\overline{\pi}}_{\gamma(0)}\equiv \SF^{\nu}(\gamma)\bmod 2$ since
$K^{\overline{\pi}}_{\gamma(0)}\equiv 0\bmod 2$. 
The  orientation for
$\mbox{\rm detind}(-K^{\overline{\pi}})=\mbox{\rm detind}(K^{\overline{\pi}})$
on the other hand is reversed by the parity of
$\dimn\kernel K^{\overline{\pi}}_{\gamma(0)}
=3$ as it's overall orientation is fixed by that at $[\mTheta]$.
The Chern--Simons functional as well as APS spectral invariants depend 
on the orientation of $Y$. Thus 
$\epsilon^{-}[A]\kappa^{-}[A]=\epsilon[A]\kappa[A]$,
$c^{-}(g,\pi')=
-c(g,\pi)$ and $\sum_{\M^{r*,-}_{\pi'}}\epsilon^{-}[A]= 
-\sum_{\M^{r*}_{\pi}}\epsilon[A]$.

Combining all of the above we obtain
$$
\tau(Y)-\tau(-Y) = c(g,\pi) +\frac{1}{2}\sum_{\M^{r*}_{\overline{\pi}}}\epsilon[A].
$$
Let $\lambda(Y)$ denote Casson's invariant \cite{AM}.
In \cite{taubes} it is established that 
$$
\frac{1}{2}\sum_{\M^{r*}_{\overline{\pi}}}\epsilon[A]=-\lambda(Y)
$$
and it was proven by Casson that $\lambda(Y)\equiv\mu(Y)\bmod 2$.
Since $c(g,\pi)\equiv\mu(Y)\bmod 2$ we obtain $\tau(Y)-\tau(-Y)\equiv 0\bmod 2$.
This completes the proof of Theorem~\ref{thrm1}.

\section{Proof of Theorem~\protect\ref{thrm2}}\label{proof2-sect}

Let us begin by reviewing the $SU(3)$--Casson invariant (in our terminology).
For more details refer to \cite{BH}. 
Denote by $\M^{SU(3)}$ the moduli space of flat $SU(3)$--connections on the trivial $SU(3)$
principal bundle over $Y$. As always $Y$ is oriented. The reducible subspace
is exactly $\M^{SU(2)}$, the moduli space of flat $SU(2)$--connections. This
coincides with $\M^{r}$ in our SW--context. A suitable class of \lq holonomy\rq\ 
perturbations $h$ can be constructed so that the perturbed space
$\M^{SU(3)}_{h}$ is non-degenerate. This means that it is a finite number
of points. Additionally each irreducible point $[A]$ has an oriented 
$\widehat{\epsilon}[A]\in\pm 1$ given by spectral flow.
However the perturbed reducible portion
$\M^{SU(3),r}_{h}$ does not consist of $SU(2)$--connections but essentially
$U(2)$--connections. $\M^{SU(3),r}_{h}$ lies in $\B_{U(2)}$ the quotient
space of $U(2)$--connections; as before there is a Chern--Simons function
$\cs$ on connnections which descends to
$\overline{\cs}\colon \B_{U(2)}\to\R/\Z$. 
To make an invariant out of $\sum_{\M^{SU(3)*}_{h}}\widehat{\epsilon}[A]$ there are
counter-terms associated to $\M^{SU(3)r*}_{h}=\M^{SU(3),r}_{h}-\{[\mTheta]\}$. However
we need to make $h$ \textit{small} which is the same condition used in our
SW--context (and from which our definition originated). Denote by $\{\N^{j}_{SU(3)}\}$ 
the corresponding system of neighbourhoods of components of 
$\M^{SU(2)}-\{[\mTheta]\}$ in $\B_{U(2)}$. 

Along the reducible
strata $\B_{U(2)}$ we have \textit{tangential} and \textit{normal}
deformation operators giving rise to tangential and normal spectral
flow quantities $\SF^{\tau}_{SU(3)}(\gamma)$ (real spectral flow), 
$\SF^{\nu}_{SU(3)}(\gamma)$ (complex spectral flow) along $\gamma$, respectively. The
term $\SF^{\tau}_{SU(3)}(\gamma)$ is used to define Taubes' orientation
$\epsilon[A]=\pm 1$ for $[A]\in\M^{SU(3)r*}_{h}$. $\SF^{\nu}_{SU(3)}(\gamma)$
is used in the term
$$
\kappa^{SU(3)}[A] = \SF^{\nu}_{SU(3)}(\gamma) +2\cs(\overline{\gamma}(0))
-2\cs(\overline{\gamma}(1)).
$$
As before $[\gamma(t)]$, $0\leq t\leq 1$ is a path from $[\mTheta]$ to 
$[A]\in \N^{j}_{SU(3)}$ say and
$[\overline{\gamma}]$ is the path from $[\mTheta]$ to the component
$\K^{j}\subset \M^{SU(2)}$, and homotopic to $[\gamma]$ rel
$\N^{j}_{SU(3)}\cup \{[\mTheta]\}$. The value of $\kappa^{SU(3)}[A]$ does
not depend on the choice of $[\gamma]$ or $[\overline{\gamma}]$.
The $SU(3)$--Casson invariant is then defined as
$$
\lambda_{SU(3)}(Y) = \sum_{\M^{SU(3)*}_{h}}\widehat{\epsilon}[A]
-\sum_{\M^{SU(3)r*}_{h}}\epsilon[A](\kappa^{SU(3)}[A]+1) \in \R.
$$
Fix a component $\K^{j}$ and homotopy class $[\gamma_{j}]$ rel $\N^{j}_{SU(3)}\cup 
\{[\mTheta]\}$ 
of paths from $[\mTheta]$ to $\N^{j}_{SU(3)}$. For every $[A]\in\N^{j}_{SU(3)}$
define $\kappa^{SU(3)}[A]$ using a path $[\gamma]$ homotopic to $[\gamma_{j}]$.
Then the Chern--Simons term is the same constant over all $[A]\in\N^{j}_{SU(3)}$, and the spectral
flow term is well-defined (depending only on $[\gamma_{j}]$). We express this
as
$$
\kappa^{SU(3)}[A] = \SF^{\nu}_{SU(3)}[A] + 2\Delta\cs(j).
$$
Thus we may rewrite the counter-term
\begin{eqnarray}
\sum_{\M^{SU(3)r*}_{h}}\epsilon[A]\kappa^{SU(3)}[A]
&=&\sum_{j} \sum_{\M^{SU(3)r*}_{h}\cap\N^{j}_{SU(3)}} \epsilon[A]\,\SF^{\nu}_{SU(3)}[A]
\label{count-equ1}\\
&&\mbox{} +2\sum_{j} \Bigl(\sum_{\M^{SU(3)r*}_{h}\cap\N^{j}_{SU(3)}}
\epsilon[A]\Bigr)\Delta\cs(j).\nonumber
\end{eqnarray}
The \textrm{local index} term
$$
\iota_{U(2)}(\K^{j}) = \sum_{\M^{SU(3)r*}_{h}\cap\N^{j}_{SU(3)}}\epsilon[A]
$$
is well-defined independent of small perturbation $h$. Given any other small non-degenerate
perturbation $h'$ we have a parameterized moduli space which is a compact oriented
cobordism between $\M^{SU(3)r*}_{h}\cap\N^{j}_{SU(3)}$ and 
$\M^{SU(3)r*}_{h'}\cap\N^{j}_{SU(3)}$.

In our SW--context make the same construction. We can identify homotopy classes
$[\gamma_{j}]$ in our SW--context with those in the $SU(3)$--Casson by the
inclusion $\B_{A}\subset\B_{U(2)}$ which is a homotopy equivalence. Then we have
in a similiar manner
\begin{eqnarray}
\sum_{\M^{r*}_{\overline{\pi}}} \epsilon[A]\,\kappa[A]
&=& \sum_{j} \sum_{\M^{r*}_{\pi}\cap\N^{j}} \epsilon[A]\,\SF^{\nu}[A]\label{count-equ2}\\
&&\mbox{} +\sum_{j} \Bigl(\sum_{\M^{r*}_{\pi}\cap\N^{j}}
\epsilon[A]\Bigr)\Delta\cs(j) \nonumber
\end{eqnarray}
and a local index
$$
\iota_{SU(2)}(\K^{j}) = \sum_{\M^{r*}_{\pi}\cap\N^{j}}\epsilon[A].
$$
The two indices $\iota_{U(2)}$ and $\iota_{SU(2)}$ are equal. This is established
by working with a restricted class of holonomy perturbations $h'$ as in
\cite{taubes} or \cite{BH} which keeps
$\M^{SU(3)r*}_{h'}$ within $SU(2)$--connections. Then it is straightforward to
relate this to our space $\M^{r*}_{\overline{\pi}}$ by a compact oriented cobordism. 
The non-integral terms for $\lambda_{SU(3)}(Y)$, $2\tau(Y)$
come from (\ref{count-equ1}), (\ref{count-equ2}) respectively.
It follows that 
$\lambda_{SU(3)}(Y)+2\tau(Y) \bmod 4$ is integral. It is also independent of 
the orientation of $Y$,
since $\lambda_{SU(3)}(Y)$ and $\tau(Y)$ are both independent of 
orientation.

\end{document}